\documentclass[11pt,amsfonts]{amsart}

\usepackage{fullpage}

\def\endproof{\qed \medskip}
\def\blacksquare{\hbox to .60em {\vrule width .60em height .60em}}

\newtheorem{theorem}{Theorem}[section]

\newtheorem{corollary}[theorem]{Corollary}

\newtheorem{lemma}[theorem]{Lemma}

\newtheorem{proposition}[theorem]{Proposition}

\newtheorem{remark}[theorem]{Remark}

\begin{document}

\title[]{On the Structure of Conformally Compact 
Einstein Metrics}

\author[]{Michael T. Anderson}

\thanks{Partially supported by NSF Grant DMS 0305865 and 0604735}

\maketitle

\abstract 
Let $M$ be an $(n+1)$-dimensional manifold with non-empty boundary, satisfying 
$\pi_{1}(M, \partial M) = 0$. The main result of this paper is that the space of 
conformally compact Einstein metrics on $M$ is a smooth, infinite dimensional 
Banach manifold, provided it is non-empty. We also prove full boundary regularity 
for such metrics in dimension 4 and a local existence and uniqueness theorem for 
such metrics with prescribed metric and stress-energy tensor at conformal infinity, 
again in dimension 4. This result also holds for Lorentzian-Einstein metrics with 
a positive cosmological constant. 
\endabstract

\section{Introduction.}
\setcounter{equation}{0}

 Let $M$ be the interior of a compact $(n+1)$-dimensional manifold $\bar M$ 
with non-empty boundary $\partial M$. A complete metric $g$ on $M$ 
is $C^{m,\alpha}$ conformally compact if there is a defining function 
$\rho$ on $\bar M$ such that the conformally equivalent metric
\begin{equation} \label{e1.1}
\widetilde g = \rho^{2}g 
\end{equation}
extends to a $C^{m,\alpha}$ Riemannian metric on the compactification 
$\bar M$. A defining function $\rho$ is a smooth, non-negative 
function on $\bar M$ with $\rho^{-1}(0) = \partial M$ and $d\rho  \neq 
 0$ on $\partial M$.

 The induced Riemannian metric $\gamma  = \widetilde g|_{\partial M}$ 
is called the boundary metric associated to the compactification 
$\widetilde g$. Since there are many possible defining functions, and 
hence many conformal compactifications of a given metric $g$, only the 
conformal class $[\gamma]$ of $\gamma$ on $\partial M$ is uniquely 
determined by $(M, g)$; the class $[\gamma]$ is called the conformal 
infinity of $g$.  Any manifold $M$ carries many conformally compact 
metrics and in this paper we are interested in Einstein metrics $g$, for 
which
\begin{equation} \label{e1.2}
Ric_{g} = -ng. 
\end{equation}
Conformally compact Einstein metrics are also called asymptotically hyperbolic 
(AH), in that $|K_{g}+1| = O(\rho^{2})$, where $K_{g}$ denotes any sectional 
curvature of $(M, g)$, at least when $g$ has a $C^{2}$ conformal compactification. 

 In this paper, we prove several distinct results on conformally 
compact Einstein metrics. First, we prove boundary regularity for such metrics 
in dimension $n+1= 4$. Thus, suppose $g$ is an Einstein metric on a 4-manifold $M$, 
which admits an $L^{2,p}$ conformal compactification, for some $p > 4$. 
If the resulting boundary metric is $C^{m,\alpha},$ or $C^{\infty}$, 
or $C^{\omega}$ (real-analytic), then $g$ is $C^{m,\alpha}$, or 
$C^{\infty}$, or $C^{\omega}$ conformally compact respectively; see 
Theorem 2.3.

 The proof of boundary regularity uses the fact that Einstein metrics 
on 4-manifolds satisfy a conformally invariant $4^{\rm th}$ order equation, 
the Bach equation, given by
\begin{equation} \label{e1.3}
\delta d(Ric - \tfrac{s}{6}g) + W(Ric) = 0. 
\end{equation}
Here $Ric$ and $g$ are viewed as 1-forms with values in $TM$, $s = tr Ric$ 
is the scalar curvature, $W$ is the Weyl curvature and $d = d^{\nabla}$ is the 
exterior derivative with $\delta  = \delta^{\nabla}$ its $L^{2}$ adjoint. The 
equations (1.3) are the Euler-Lagrange equations for the $L^{2}$ norm of the 
Weyl curvature $W$, as a functional on the space of metrics on $M$. Since 
(1.3) is conformally invariant in dimension 4, any conformal compactification 
$\widetilde g$ of $g$ satisfies (1.3) and boundary regularity is established 
by studying the boundary regularity of solutions of the non-degenerate equation 
(1.3) on $\bar M$. Theorem 2.3 corrects a small gap in the proof of boundary 
regularity in [2, Thm.2.4], cf.~Remark 2.4. 

 As an application of these techniques, we also prove a local existence and 
uniqueness result. Thus, recall the Fefferman-Graham expansion of an AH Einstein 
metric [12]; in dimension 4, this is given by
\begin{equation} \label{e1.4}
\bar g = t^{2}g \sim  dt^{2} + g_{(0)} + t^{2}g_{(2)} + t^{3}g_{(3)} + 
\cdots + t^{k}g_{(k)} + ... , 
\end{equation}
where $t$ is a geodesic defining function, i.e. $t(x) = dist_{\bar 
g}(x, \partial M)$. The boundary metric $\gamma$ is given by $\gamma  
= g_{(0)}$ and the term $g_{(2)}$ is intrinsically determined by 
$\gamma$. The Einstein constraint equations at conformal infinity $\partial M$ 
are equivalent to the statement that the term $g_{(3)}$ is transverse-traceless 
on $(\partial M, \gamma)$, i.e.~$\delta_{\gamma}g_{(3)} = tr_{\gamma}g_{(3)} = 0$, 
see for instance [11]. However, beyond this, the $g_{(3)}$ term is not determined 
by the boundary metric $\gamma$. All higher order terms in the expansion (1.4) are 
determined by $g_{(0)}$ and $g_{(3)}$ via the Einstein equations. It is also 
worth noting that from a physics perspective, the term $g_{(3)}$ is identified 
with the stress-energy tensor of the conformal infinity, cf.~again [11] for 
instance.

 In [12], Fefferman-Graham proved that if $\gamma  = g_{(0)}$ is any 
real-analytic metric on an arbitrary 3-manifold $\partial M$, and one 
sets $g_{(3)} = 0$, so that the formal expansion (1.4) is even in $t$, 
then there exists a real-analytic AH Einstein metric defined in a 
thickening $M = [0,\varepsilon )\times \partial M$, with boundary metric 
$\gamma$. Thus, the series (1.4) converges to $\bar g$. This was proved by 
using results of Baouendi-Goulaouic on the convergence of formal series 
solutions to nonlinear Fuchsian systems of PDE's. A result analogous to 
this was proved earlier by LeBrun [19] for self-dual Einstein metrics 
on thickenings of 3-manifold boundaries, using twistor methods. 

 The following result generalizes these results to allow for an 
arbitrary $g_{(3)}$ term.
\begin{theorem} \label{t 1.1.}
  Let $N$ be a closed 3-manifold, and let $(\gamma, \sigma)$ be a 
pair consisting of a real-analytic Riemannian metric $\gamma$ on $N$, 
and a real-analytic symmetric bilinear form $\sigma$ on $N$ satisfying 
$\delta_{\gamma}\sigma = tr_{\gamma}\sigma  = 0$. Then there exists a 
unique (up to isometry), $C^{\omega}$ conformally compact Einstein metric 
$g$, defined on a thickening $N\times I$ of $N$, for which the expansion 
\eqref{e1.4} converges to $\bar g$ and is given by
\begin{equation} \label{e1.5}
\bar g = dt^{2} + \gamma  + t^{2}g_{(2)} + t^{3}\sigma  + ... + 
t^{k}g_{(k)} + ... 
\end{equation}
\end{theorem}

 The proof is based again on the Bach equation, together with the 
Cauchy-Kovalewsky theorem. An analogous result also holds for 
Lorentzian-Einstein metrics, i.e.~solutions to the vacuum Einstein 
equations in general relativity with a positive cosmological constant 
$\Lambda$, cf.~Theorem 2.6. The result in this case is related to work 
of H. Friedrich [13].

\medskip

 Next, we turn to the structure of the moduli space of AH Einstein 
metrics on a given $(n+1)$-manifold $M$. Let $E_{AH} = E_{AH}^{m,\alpha}$ 
be the space of AH Einstein metrics $g$ on $M$ which admit a $C^{m,\alpha}$ 
compactification $\widetilde g$ as in (1.1). We require that $m \geq 2$, 
$\alpha\in (0,1)$ but otherwise allow any value of $m$, including $m = 
\infty$ or $m = \omega$. The space $E_{AH}^{m,\alpha}$ is given the 
$C^{m,\mu}$ topology on $\bar M$, for any fixed $\mu < \alpha$, via a fixed 
compactification as in (1.1). Let ${\mathcal E}_{AH} = E_{AH}/${\rm Diff}$_{1}(\bar M)$, 
where ${\rm Diff}_{1}(\bar M) = {\rm Diff}_{1}^{m+1,\alpha}(\bar M)$ is the 
group of $C^{m+1,\alpha}$ diffeomorphisms of $\bar M$ inducing the identity 
on $\partial M$, acting on $E_{AH}$ in the usual way by pullback. 

  Regarding the boundary data, let $Met(\partial M) = Met^{m,\alpha}(\partial M)$ 
be the space of $C^{m,\alpha}$ metrics on $\partial M$ and ${\mathcal C}  = {\mathcal C} 
(\partial M)$ the corresponding space of pointwise conformal classes, 
endowed with the $C^{m,\mu}$ topology as above. There is a natural boundary 
map, (for any fixed $(m, \alpha)$),
\begin{equation} \label{e1.6}
\Pi : {\mathcal E}_{AH} \rightarrow  {\mathcal C} , \ \ \Pi [g] = [\gamma], 
\end{equation}
which takes an AH Einstein metric $g$ on $M$ to its conformal infinity 
on $\partial M$.

 We then have the following result on the structure of ${\mathcal E}_{AH}$ 
and the map $\Pi$.
\begin{theorem} \label{t 1.2.}
  Let $M$ be a compact, oriented 4-manifold with boundary $\partial M$ 
satisfying $\pi_{1}(M, \partial M) = 0$. If, for a given $(m,\alpha)$, 
$m \geq 3$, ${\mathcal E}_{AH}$ is non-empty, then ${\mathcal E}_{AH}$ 
is a $C^{\infty}$ smooth infinite dimensional separable Banach manifold. 
Further, the boundary map 
\begin{equation} \label{e1.7}
\Pi : {\mathcal E}_{AH} \rightarrow  {\mathcal C}  
\end{equation}
is a $C^{\infty}$ smooth Fredholm map of Fredholm index 0.
\end{theorem}

 Implicit in Theorem 1.2 is the boundary regularity statement that an 
AH Einstein metric with $C^{m,\alpha}$ conformal infinity has a 
$C^{m,\alpha}$ compactification. Versions of Theorem 1.2 also hold in 
arbitrary dimensions $n > 4$; see Theorems 5.5 and 5.6 for the 
precise statements. 

  The condition $\pi_{1}(M, \partial M) = 0$ is equivalent to the statements 
that $\partial M$ is connected and the inclusion map $\iota: \partial M 
\rightarrow M$ induces a surjection $\pi_{1}(\partial M) \rightarrow \pi_{1}(M) 
\rightarrow 0$. It is not clear whether Theorem 1.2 holds globally without 
this assumption, although generic metrics in ${\mathcal E}_{AH}$ always have 
smooth manifold neighborhoods, cf.~Remark 3.2 and the discussion following 
Theorem 4.1. 

 Theorem 1.2 is a generalization of previous results of Graham-Lee [15] 
and Biquard [8], who proved a local analogue of this result, (without 
full boundary regularity), in neighborhoods of metrics $g\in{\mathcal E}_{AH}$ 
which are regular points of the map $\Pi$. The proof of Theorem 1.2 uses methods 
introduced in [15] and [8]. Further, Theorem 1.2 is formally analogous to results 
on the space of minimal surfaces, cf.~[9] and especially [26], [27] and we have 
also been influenced by this work.

  The proof of Theorem 1.2 requires a rather subtle understanding of the behavior 
of infinitesimal AH Einstein deformations which are in $L^{2}(M, g)$; one needs to 
know that such deformations satisfy a suitable unique continuation property at 
infinity. This was proved in [6] and is presented here in \S 3, cf.~Proposition 3.1, 
after some preliminary introductory material. The main work in the proof of 
Theorem 1.2 is then given in \S 4, with the final proof given in \S 5. 
We also point out that it is proved in Theorem 5.7 that the spaces 
${\mathcal E}_{AH}^{m,\alpha}$ are stable in $(m, \alpha)$; they are all 
diffeomorphic and the inclusion of ${\mathcal E}_{AH}^{m',\alpha'}$ into 
${\mathcal E}_{AH}^{m,\alpha}$ for $m' +\alpha'  > m+\alpha $ is dense, 
including the case $m' = \infty$ or $m' = \omega$. Thus, the structure of 
the spaces ${\mathcal E}_{AH}^{m,\alpha}$ is essentially independent of 
$(m, \alpha)$.

\medskip

   The results in this paper are also used in [5], which studies the existence 
problem for conformally compact Einstein metrics with prescribed conformal 
infinity on 4-manifolds.

\section{The Bach equation and AH Einstein metrics.}
\setcounter{equation}{0}

 In this section, we prove boundary regularity for AH Einstein metrics 
in dimension 4, together with Theorem 1.1 and various applications.

\medskip

 We begin with the study of boundary regularity. Let $g$ be an AH 
Einstein metric on a 4-manifold $M$. Then $g$ satisfies the conformally 
invariant Bach equation (1.3). Hence, any conformal compactification 
$\widetilde g$ of $g$ also satisfies (1.3). In the following, to 
simplify notation, we work with a given conformal compactification 
$\widetilde g$ of an AH Einstein metric $g$ and drop the tilde from 
the notation; thus, from now on until Corollary 2.2, $g$ denotes $\widetilde g$. 

 By a standard Weitzenbock formula, (1.3) may be rewritten in the form
\begin{equation} \label{e2.1}
D^{*}DRic = -\tfrac{1}{3}D^{2}s - \tfrac{1}{6}\Delta s + {\mathcal R} , 
\end{equation}
where ${\mathcal R}$ is a term quadratic in the curvature of $g$.

 As it stands, the equation (2.1) (or (1.3)) does not form an elliptic 
system, due to its invariance under diffeomorphisms and conformal 
deformations. Since we wish to cast (2.1) in the form of an elliptic 
boundary value problem, two choices of gauge are needed to break these 
symmetries.

 First, with regard to the diffeomorphism invariance, we use, as is now 
common, harmonic coordinates. Thus let $x^{i}$, $i = 1,2,3$ be local 
harmonic coordinates on $(\partial M, \gamma)$ and extend $x^{i}$ 
locally into $M$ by requiring that $x^{i}$ is harmonic with respect to 
$g$ (i.e.~$\widetilde g$):
\begin{equation} \label{e2.2}
\Delta x^{i} = 0, 
\end{equation}
where the Laplacian is with respect to $(M, g)$. Also, let $x^{0}$ be a local 
``harmonic defining function'', satisfying
\begin{equation} \label{e2.3}
\Delta x^{0} = 0, \ \  x^{0}|_{\partial M} = 0 .
\end{equation}
Thus, the functions $x^{i}$, $i = 0,1,2,3$, form a local coordinate 
system for $M$ up to its boundary, harmonic with respect to $g$.

  Suppose in a given fixed (or background) $C^{\infty}$ atlas for $M$ near 
$\partial M$, the compactified metric $g$ is in $L^{k,p}$, for some $k \geq 2$, 
$p > 4$, or in $C^{m,\alpha}$, $m \geq 2$. Then it is well-known, and will be 
frequently be used below, that $g$ is $L^{k,p}$ or $C^{m,\alpha}$ respectively 
in local boundary $g$-harmonic coordinates. Further, these harmonic coordinates 
are $L^{k+1, p}$ or $C^{m+1,\alpha}$ functions of the background local coordinates 
respectively. The same remarks pertain in the real-analytic case. Thus, harmonic 
coordinates give optimal regularity properties. 

  With regard to the conformal invariance, it is natural to specify the 
scalar curvature $s$ to determine the conformal gauge. At a later 
point, we will choose a Yamabe gauge, where $s = const$. However, for 
the moment, we assume that $s$, the scalar curvature of $g$, is a given 
function with a given degree of smoothness. In particular, the two 
scalar curvature terms on the right in (2.1) are thus ``determined''. 
In the following, Greek indices $\alpha$,$\beta$ run over $0,1,2,3$ while 
Latin indices $i$, $j$ run over $1,2,3$.

 In local harmonic coordinates, $-2D^{*}DRic = \Delta\Delta 
g_{\alpha\beta}$ + ($3^{\rm rd}$ order terms). Thus, the system (2.1) 
may be rewritten in local harmonic coordinates as
\begin{equation} \label{e2.4}
\Delta\Delta g_{\alpha\beta} = F_{\alpha\beta}(g, s), 
\end{equation}
where $F$ is of order 3 in $g$, order 2 in $s$, with real-analytic 
coefficients; here $\Delta  = 
g^{\alpha\beta}\partial_{\alpha}\partial_{\beta}$ and $s$ is treated as given. 

 This is a $4^{\rm th}$ order elliptic system, with leading order term in 
diagonal form.

 We now set up the boundary conditions for this system. On the surface, 
it would be simplest to just choose Dirichlet and Neumann boundary 
conditions on the set of all $g_{\alpha\beta}$. However, via the map $\Pi$ in 
(1.6), at the boundary we only have information on the intrinsic metric 
$\gamma_{ij} = g_{ij}$; the $g_{0\alpha}$ terms are gauge dependent, and have 
no apriori prescribed form at $\partial M$. Thus, it is not clear if 
$\{g_{0\alpha}\}$ satisfy any particular boundary conditions apriori. 
Moreover, in general Bach-flat metrics are not conformally Einstein; conformally 
Einstein metrics thus necessarily induce certain special boundary conditions. In 
the case of geodesic gauge, this is discussed in the proof of Theorem 1.1 below, 
and is closely related to the Fefferman-Graham expansion. However, such 
a gauge is badly behaved for elliptic boundary value problems. While it 
is an interesting open question to characterize the boundary conditions 
for Bach-flat metrics to be conformally Einstein (in arbitrary gauges), this 
issue will not be addressed here; instead we will derive certain boundary 
conditions for the Bach-flat equations as a consequence of the metric being 
conformally Einstein. 

 To begin, we divide the collection $\{g_{\alpha\beta}\}$ into two parts. 
First $g_{ij}$ are the tangential components of $g$, with $1 \leq i,j \leq 3$. 
The remaining terms $g_{0\alpha}$ are the mixed and normal terms. The basic idea 
is then to impose Dirichlet boundary conditions on $g_{ij}$, of $0^{\rm th}$ and 
$2^{\rm nd}$ order, while imposing Neumann-type boundary conditions on $g_{0\alpha}$, 
of $1^{\rm st}$ and $3^{\rm rd}$ order. These latter two conditions come from the 
gauge choice of harmonic coordinates. 

 In more detail:

 $\{1\}$. Dirichlet boundary conditions on $g_{ij}$:
\begin{equation} \label{e2.5}
g_{ij} = \gamma_{ij} \ \ {\rm on} \ \  \partial M \ \ {\rm (locally)}, 
\end{equation}
where $\gamma $ is the given boundary metric on $\partial M$.

 $\{2\}$. Neumann-type boundary conditions on $g_{0\alpha}$. It is 
convenient to set these up for the inverse variables $g^{0\alpha}$. 
These are of the form
\begin{equation} \label{e2.6}
N(g^{00}) = -2Hg^{00}, 
\end{equation}
\begin{equation} \label{e2.7}
N(g^{0i}) = \tfrac{1}{2}q^{i\beta}\partial_{\beta}g^{00} 
- Hg^{0i}. 
\end{equation}
Here $N = (g^{00})^{-1/2}g^{0\beta}\partial_{\beta} = q^{0\beta}\partial_{\beta}$ 
is the unit normal, and $q^{\alpha\beta} = (g^{00})^{-1/2}g^{\alpha\beta}$. The term 
$H$ is the mean curvature, $H = g^{ij}A_{ij}$, with $A_{ij} = 
\frac{1}{2}q^{0\alpha}[\partial_{\alpha}g_{ij} - (\partial_{i}g_{\alpha j} + 
\partial_{j}g_{\alpha i})]$ the $2^{\rm nd}$ fundamental form of the boundary. 
The fact that the components $\{g^{0\alpha}\}$ satisfy (2.6)-(2.7) in boundary 
harmonic coordinates $\{x_{\alpha}\}$ was derived in [2], cf.~also [3]. This 
does not require the metric to be conformally Einstein or Bach-flat; it holds 
in general.

  The equations (2.6)-(2.7) can be reexpressed as Neumann-type conditions on 
the coefficients $g_{0\alpha}$, since $g_{ij}$ is given on $\partial M$ by 
(2.5). However, these expressions will not be given explicitly, since only the 
linearized versions of (2.6)-(2.7) need to be actually computed. 

 $\{3\}$. Dirichlet boundary conditions on $Ric_{ij}$. One has, in harmonic 
coordinates in general, $-2Ric_{ij} = \Delta^{M}g_{ij} + Q(g, \partial g)$, 
where the Laplacian is with respect to $(M, g)$ and $Q$ is a $1^{\rm st}$ order term 
in $g$. A simple and standard calculation shows that for conformally Einstein metrics 
$g$, $Ric_{ij}$ is determined by $(Ric_{\gamma})_{ij}$ at $\partial M$, modulo 
lower order terms. In fact, cf.~ [2,Lemma 1.3] for instance, the extrinsic and 
intrinsic Ricci curvatures are related by
\begin{equation} \label{e2.8}
Ric_{ij} = 2(Ric_{\gamma})_{ij} + \tfrac{1}{6}(s - 
\tfrac{3}{2}s_{\gamma})\gamma_{ij} - (\tfrac{H}{3})^{2}\gamma_{ij}. 
\end{equation}
Since also $-2(Ric_{\gamma})_{ij} = \Delta^{\partial M}\gamma_{ij} + Q(\gamma , 
\partial\gamma)$, this gives
\begin{equation} \label{e2.9}
\Delta^{M}g_{ij} = F(\gamma , \partial\gamma ,\partial^{2}\gamma , H, s) 
\ \ {\rm on} \ \ \partial M, 
\end{equation}
where $F$ is real-analytic in its arguments; observe that $H$ is 
real-analytic in $g$ and its $1^{\rm st}$ derivatives.

 $\{4\}$. Neumann-type boundary conditions on $Ric_{0\alpha}$. These 
turn out to be
\begin{equation} \label{e2.10}
q^{0\beta}N(\Delta^{M}g_{0\beta}) + {\tfrac{1}{3}}\Delta^{\partial M}H = 
-\partial_{0}s + Q_{0}, 
\end{equation}
\begin{equation} \label{e2.11}
q^{0\beta}N(\Delta^{M}g_{i\beta}) = -{\tfrac{2}{3}}\partial_{i}s + Q_{i}, 
\end{equation}
where $Q = (Q_{0}, Q_{i})$ is an operator of order less than 3. Recall 
that the scalar curvature $s$ is treated as given. 

   The equations (2.10)-(2.11) are basically a consequence of the Bianchi identity. 
To derive them, the contracted Bianchi identity gives $\delta Ric = -\frac{1}{2}ds$, 
so that for any $\alpha$, $-\delta Ric(\partial_{\alpha}) = 
\frac{1}{2}\partial_{\alpha}s$. By definition, $-\delta Ric(\partial_{\alpha}) = 
(\nabla_{N}Ric)(N,\partial_{\alpha}) + (\nabla_{e_{j}}Ric)(e_{j},\partial_{\alpha})$, 
where $e_{j}$ runs over an orthonormal basis tangent to $\partial M$. Also 
$(\nabla_{N}Ric)(N,\partial_{\alpha}) = N(Ric(N,\partial_{\alpha})) - 
Ric(N, \nabla_{N}\partial_{\alpha}) - Ric(\nabla_{N}N, \partial_{\alpha})$, 
and $N(Ric(N,\partial_{\alpha})) = q^{0\beta}N(Ric_{\alpha\beta}) + 
N(q^{0\beta})Ric_{\alpha\beta}$. Putting these together, the Bianchi identity 
may be rewritten as
\begin{equation} \label{e2.12}
q^{0\beta}N(Ric_{\alpha\beta}) + e_{j}(Ric(e_{j},\partial_{\alpha})) 
=\tfrac{1}{2}\partial_{\alpha}s + Q_{\alpha}  ,
\end{equation}
where $Q_{\alpha}$ involves the Ricci curvature of $g$ and the $1^{\rm st}$ derivatives 
of $g$ at $\partial M$; thus $Q_{\alpha}$ is of order less than 3.  Using again the fact 
that $-2Ric_{\alpha\beta} = \Delta^{M}g_{\alpha\beta} + Q(g, \partial g)$ in 
harmonic coordinates, the equation above may be rewritten as
\begin{equation} \label{e2.13}
q^{0\beta}N(\Delta^{M}g_{0\beta}) + div^{\partial M}(\Delta^{M}g_{0\cdot}) = 
-\partial_{0}s + Q_{0}. 
\end{equation}
\begin{equation} \label{e2.14}
q^{0\beta}N(\Delta^{M}g_{i\beta}) + div^{\partial M}(\Delta^{M}g_{i\cdot}) = 
-{\tfrac{2}{3}}\partial_{i}s + Q_{i}, 
\end{equation}
where we have separated the cases $\alpha = 0$ and $\alpha = i > 0$. 

  These equations have to be modified somewhat, since, together with 
$\{1\} - \{3\}$, they do not lead to elliptic boundary conditions; this is probably 
because both systems (2.6)-(2.7) and (2.13)-(2.14) are identities in harmonic 
coordinates. In any case, we reexpress the divergence terms in (2.13) and (2.14) 
as follows. Regarding (2.13), a standard formula for the change of Ricci curvature under 
conformal changes gives, at $\partial M$, 
\begin{equation} \label{e2.15}
Ric_{0i} = -\frac{2}{3}\partial_{i}H , 
\end{equation}
cf.~again [2,~(1.18)] for instance. This equation uses the fact that $g$ is conformally 
Einstein; note this is $2^{\rm nd}$ order in $g_{\alpha\beta}$, in contrast to (2.8), 
(cf.~also Remark 2.4). Substituting (2.15) in (2.12) and using (2.13) then gives (2.10). 
For (2.14), the ambient Ricci curvature term $Ric_{ij}$ in (2.12), ($\alpha = i > 0$), 
is in fact intrinsic to the boundary, modulo lower order terms and the $s$ term, by (2.8). 
Since the derivatives are also being taken tangentially, the term 
$e_{j}(Ric(e_{j}, \partial_{i}))$ is intrinsic to the boundary metric $\gamma$, 
modulo lower order terms and the $s$ term. At leading order, it does not depend on 
$g_{0\alpha}$ and thus it may be absorbed into the $Q_{i}$ term. Taking the $s$ term 
in (2.8) into account, this gives (2.11). 

\medskip

 The boundary conditions $\{1\}-\{4\}$ are the conditions that will be used below. 
Note that only the conditions $\{3\}-\{4\}$ use the fact that $g$ is conformally 
Einstein. Given this groundwork, we are now in position to prove the following 
result.
\begin{proposition} \label{p 2.1.}
  Let $\widetilde g$ be a conformal compactification of an AH Einstein 
metric g, defined on a thickening $M$ of $N = \partial M$, with scalar 
curvature $s = s_{\widetilde g}$ given. In boundary harmonic coordinates, 
the Bach equation (2.1), with the boundary conditions $\{1\}-\{4\}$, forms a 
non-linear elliptic boundary value problem, with real-analytic coefficients.
\end{proposition}
{\bf Proof:}
 It is clear that the operator (2.1) and boundary operators $\{1\}-\{4\}$ are 
real-analytic in the metric $g$ and its derivatives. Thus, one needs to 
check that the conditions of Agmon-Douglis-Nirenberg [1] or Morrey [23,~\S 6] 
are satisfied; we will verify the conditions of Morrey. First, 
ellipticity of the boundary value problem depends only on that of its 
linearization at any solution $\widetilde g.$ Thus, in the work above, 
replace $\widetilde g$ by $\widetilde g + \lambda \widetilde h$ and take the 
derivative with respect to $\lambda$ to obtain a linear system in $\widetilde h$; as 
above, henceforth we drop the tilde from the notation. We will also 
assume that the coordinate system is small, so that $g_{\alpha\beta}$ 
is close to $\delta_{\alpha\beta};$ in particular $g^{\alpha\beta} \sim  
g_{\alpha\beta}$.

 The interior system is then essentially the same as before: 
\begin{equation} \label{e2.16}
\Delta\Delta h + F_{3}(g, h) = 0, 
\end{equation}
where the Laplacian is with respect to $(M, g)$ and $F$ is of order 3 in $h$.

 In the notation of Morrey [23,\S 6.1], the interior system has the form
$$L_{jk}u^{k} = 0 \ \ {\rm in} \ \  M.  $$
Here $\{u^{k}\} = \{h_{\alpha\beta}\}$ so that $j, k \in  \{1, ..., N\}$, 
with $N = \#(\alpha\beta ) = 10$. The leading order term $L_{jk}'$ of 
$L_{jk}$ is given by $L_{jk}'  = (\Delta\Delta )\delta_{jk},$ the 
biLaplacian acting diagonally. 

 The order of each $L_{jk}$ is 4, and we set $t_{k} = 4$, for all $k$, 
$s_{j} = 0$, for all $j$. This leading order symbol of $L_{jk}$ has $2m$ 
roots, each $+i$ or $-i$ (at a cotangent vector of the form $\xi +n$, 
where $\xi$  is tangent to $\partial M$, $|\xi| = 1$, and $n$ is the unit 
conormal. Hence, the system (2.16) is properly elliptic, (i.e. satisfies the 
root condition).

 The boundary operator has the general form
\begin{equation} \label{e2.17}
B_{rk}u^{k} = f_{r} \ \ {\rm on} \ \ \partial M,
\end{equation}
where $r \in \{1, ..., m\}$, with $m = 2N$. Thus, one has 2 boundary 
operators for each $h_{\alpha\beta}$. The operator $B_{rk}$ is considered 
as a $2N\times N$ matrix, with each $N\times N$ block consisting of 6 
horizontal rows $(ij)$ ordered lexicographically and 4 mixed rows, 
ordered $(00, 01, 02, 03)$.

 The order of $B_{rk} = t_{k}-h_{r} = 4-h_{r}$. Thus, for $u^{k} = 
\{h_{ij}\}$, i.e. the 6 tangential components of $h$, one has $h_{r} = 4$ 
for $1 \leq  r \leq  6$, (corresponding to the Dirichlet data $\{1\}$) and 
$h_{r} = 2$ for $11 \leq  r \leq 16$, (corresponding to the Dirichlet 
data $\{3\}$). Similarly, for the mixed terms $h_{0\alpha},$ one has $h_{r} 
= 3$ for $7 \leq  r \leq 10$ and $h_{r} = 1$ for $17 \leq  r \leq  
20$, corresponding to the Neumann data $\{2\}$ and $\{4\}$. Thus, $h_{0} = 0$ in 
the notation of [23,\S 6.1].

 Since $g_{\alpha\beta} \sim  \delta_{\alpha\beta},$ the positive roots 
$z_{s}^{+}(x, \xi )$, $s = 1, ... m = 20$, of the symbol $L(x, \xi 
+zn)$ are all close to $z = i|\xi|$; when $g_{\alpha\beta} = 
\delta_{\alpha\beta},$ the positive roots are exactly $z = i|\xi|.$ Hence
$$L_{0}^{+}(x, \xi , z) \equiv \prod_{1}^{20}(z - z_{s}^{+}(x,\xi )) \sim  
(z - i|\xi|)^{20}. $$
Let $L^{jk}$ be the matrix adjoint of $L_{jk}$, (the matrix of minors of $L$). 
Then $L^{jk}(x, \xi +zn) = (|\xi|^{2}+|z|^{2})^{18}\delta_{jk}$ + lower order 
terms. One then forms
\begin{equation} \label{e2.18}
Q_{rk}(x,\xi ,z) = \sum_{j=1}^{N}B_{rj}' (x,\xi +zn)L^{jk}(x,\xi +zn) 
\sim  B_{rk}' (x,\xi +zn)(|\xi|^{2}+|z|^{2})^{18}, 
\end{equation}
where $B_{rk}'$ is the leading order symbol of $B_{rk}$. 

 Here $Q_{rk}$ is viewed as a polynomial in $z$, for any fixed $x$, $\xi$ 
with $x\in\partial M$ and $\xi$ a cotangent vector to $\partial M$; $n$ 
is the unit conormal. Then the complementary condition is that the rows 
of $Q_{rk}$ are linearly independent mod $L_{0}^{+},$ i.e. 
$$\sum_{r=1}^{m}c_{r}Q_{rk}(x,\xi ,z) = 0 \ \ mod \ \ L_{0}^{+} \Rightarrow  
\{c_{r}\} = 0. $$
By (2.18), this is essentially equivalent to 
\begin{equation} \label{e2.19}
\sum_{j=1}^{m}c_{r}B_{rk}' (x,\xi +zn) \equiv (z - i|\xi|)^{2} 
\Rightarrow  \{c_{r}\} = 0, 
\end{equation}
where the congruence is modulo polynomials in $z$. More precisely, this is the 
condition one obtains when the lower order terms in $L^{jk}$ are ignored. 
Including the lower order terms leads to the addition of polynomials of higher 
degree on the right in (2.19), and it will be obvious from the computations below 
that one may safely ignore such terms. 

   The $2N\times N$ matrix $B_{rk}'$ is the leading order symbol for the 
linearization of the boundary problems $\{1\}$-$\{4\}$ at $g$, with variable or 
unknown $h$. Consider this matrix $M$ as a pair of $N\times N$ matrices, 
an upper block $M_{1}$ consisting of the boundary operators $\{1\}$ and $\{2\}$, 
and a lower block $M_{2}$, consisting of the boundary operators $\{3\}$ and $\{4\}$. 
The leading order symbol of the linearization of $B$ is obtained by 
replacing $g$ by $h$ in the highest derivatives of $g$ that appear in 
$\{1\}$-$\{4\}$, and ignoring all lower order terms. Further, since the ellipticity 
condition is open and we are working locally, one may assume that $g_{\alpha\beta} 
= \delta_{\alpha\beta}$. A simple inspection of the form of $\{1\}$-$\{4\}$ then 
leads to the following description of $M$.

 The matrix $M_{1}$ consists of $I_{6}$, the $6\times 6$ identity matrix, with 
0 elsewhere in the first 6 rows, corresponding to the boundary operator 
$\{1\}$. For the next 4 rows, corresponding to the boundary operator $\{2\}$, 
the $4\times 4$ block corresponding to the $(0\alpha)$ terms (the lower right block) 
has the form 
\begin{equation} \label{e2.20}
\left(
\begin{array}{cccc}
z & 2\xi_{1} & 2\xi_{2} & 2\xi_{3} \\
\frac{1}{2}\xi_{1} & z & 0 & 0 \\
\frac{1}{2}\xi_{2} & 0 & z & 0 \\
\frac{1}{2}\xi_{3} & 0 & 0 & z 
\end{array}
\right)
\end{equation}
To see this, the diagonal $z$ terms come from the operator $N$ in (2.6)-(2.7). 
Next, one has $H = g^{ij}A_{ij} = A_{ii} = -\partial_{i}h_{0i} + 
\frac{1}{2}\partial_{0}h_{ii}$. Via (2.6), the first term here gives rise 
to the first row in (2.20); the term $\frac{1}{2}\partial_{0}h_{ii}$, giving rise 
to $-zh_{ii}$ goes into the $(ii)$ columns of the $(00)$ row, (in the lower left 
block), and may ignored. The first $(00)$ column in (2.20) comes from the first 
term on the right in (2.7); note that the term $Hg_{0i}$ linearizes to 0. 

   The matrix $M_{2}$ has a similar description. The first 6 rows of 
$M_{2}$ consist of $(|z|^{2}+|\xi|^{2})I_{6}$ and 0 elsewhere, 
corresponding to the boundary operator $\{3\}$. The last $4\times 4$ block (on 
the lower right) of the boundary operator $\{4\}$ gives a matrix of the form 
\begin{equation} \label{e2.21}
\left(
\begin{array}{cccc}
z(z^{2} + |\xi|^{2}) & -{\tfrac{1}{3}}|\xi|^{2}\xi_{1} & 
-{\tfrac{1}{3}}|\xi|^{2}\xi_{2} & -{\tfrac{1}{3}}|\xi|^{2}\xi_{3} \\
0 & z(z^{2} + |\xi|^{2}) & 0 & 0 \\
0 & 0 & z(z^{2} + |\xi|^{2}) & 0 \\
0 & 0 & 0 & z(z^{2} + |\xi|^{2}) 
\end{array}
\right)
\end{equation}
There are again terms of the form $|\xi|^{2}z$ in the tangential $(ii)$ columns of 
the $(00)$ row, (in the lower left block) but again these can be ignored. 

  Using these forms of the matrices, together with the fact that the terms ignored 
above are at most first order in $z$, it is easy to see that there are no non-trivial 
solutions of (2.19). Namely, the polynomials in $z$ on the left side of (2.19) are 
all of order at most 3, with no $z^{2}$ terms. Such polynomials cannot have 
$i|\xi|$ as a double root. This shows that all the hypotheses of [23,~\S6.1,6.3] 
are satisfied, which proves the result.

{\endproof}

 Having verified that the Bach equation with the boundary conditions 
$\{1\}$-$\{4\}$ forms an elliptic boundary value problem in local harmonic 
coordinates, one then has the following:

\begin{corollary} \label{c 2.2.}
  Let $g$ be an AH Einstein metric on a 4-manifold $M$, which admits an 
$L^{2,p}$ conformal compactification $\widetilde g$ in local boundary 
harmonic coordinates, for some $p > 4$, with boundary metric $\gamma$.

 Let $k \geq 1$ and $q \geq 2$. If $\gamma\in L^{k+2,q}(\partial M)$ and 
the scalar curvature $\widetilde s \in  L^{k,q}(M)$, with 
$\widetilde s|_{\partial M} \in  L^{k,q}(\partial M)$, then the metric 
$\widetilde g \in  L^{k+2,q}(M)$. 

 Similarly, for $m \geq 0$ and $\alpha \in (0,1)$, if $\gamma \in 
C^{m+2,\alpha}(\partial M)$ and $\widetilde s \in C^{m,\alpha}(M)$, 
then $\widetilde g \in C^{m+2,\alpha}(M)$. If $\gamma$ and $\widetilde s$ 
are real-analytic, then so is $\widetilde g$. 
\end{corollary}
{\bf Proof:} 
 The regularity hypotheses and conclusions are understood to be with respect 
to local boundary harmonic coordinates. 

 This result follows from Proposition 2.1 and the regularity theory for elliptic 
systems, cf.~[23, \S 6]. Suppose first that the compactification $\widetilde g$ is 
$L^{4,p}$ or $C^{4,\alpha}$, so that $\widetilde g$ is a classical 
solution of the Bach equation (2.1) and satisfies the boundary conditions 
$\{1\}$-$\{4\}$, with the given control on $\widetilde s$. Then the result 
follows from boundary regularity for such elliptic systems, see 
[23, Thm.~6.3.7]. Here, the coefficients of the interior operator $L$ and 
boundary operator $B$ are frozen to obtain a linear elliptic system, and 
the usual boostrap argument is used to obtain regularity. In the notation of 
[23], one sets $h_{0} = 0$, $h = 1$, and proceeds iteratively. The real-analytic 
case follows from [23, Thm.~6.7.6$'$]. 

 In the case $\widetilde g \in L^{2,p}$, the metric $\widetilde g$ 
satisfies the Bach equation and boundary conditions weakly. Using the fact that 
these equations are of divergence-form at leading order, (since they come from 
a natural variational problem), one then applies [23, Thm.~6.4.8] to prove that 
$\widetilde g$ has higher regularity, according to the regularity of $\widetilde s$ 
and the boundary data. This process is then iterated until $\widetilde g$ is a 
classical solution, as above. 

  To verify this in more detail, in the Bach equation (2.1), the assumptions 
$s \in L^{2,q}$ and $g \in L^{2,p}$ imply that the right side of (2.1) is in 
$L^{\hat q}$, where $\hat q = q/2 > 1$. Hence, 
$$D^{*}DRic = f \in L^{\hat q}.$$ 
(In Morrey's notation, $\hat q$ equals $q$ and $f$ is a term $f_{j}^{0}$ in 
[23, (6.4.1)]). Next, the leading term of $D^{*}DRic$ is the biLaplacian 
$\Delta \Delta$. In local harmonic coordinates, the expression $(\Delta \Delta u)dV$ 
schematically has the form $\partial(g \partial(g\partial^{2}u)) = \partial^{2}
(g^{2}\partial^{2}u) - \partial(g\partial g \partial^{2}u)$; here $g$ or $g^{2}$ 
denotes some algebraic expression in the metric $g$. The interior operator $L(u)$ 
is now chosen to be the first term, (the leading order term of the biLaplacian),
$$L(u) = \partial^{2}(g^{2}\partial^{2}u),$$
while the second term $\partial(g\partial g \partial^{2}u)$ is treated as a term 
$f_{j}^{1}$ in [23, (6.4.1)]. Note that $g\partial g \partial^{2}u \in L^{\hat q}$ 
again, when one sets $u = g$. The remaining lower order terms in $D^{*}DRic$ are then 
treated in exactly the same way; it should be noted here that all $3^{\rm rd}$ order 
terms in the metric in $D^{*}DRic$ are total derivatives, i.e.~of divergence form, 
as above. 

  In Morrey's notation, one now chooses, $h = h_{0} = -2$, and $h' = -1$ and for 
the interior operator sets, $m_{j} = 2$ and $s_{j} = 0$ for all $j$, while $t_{k} 
= 4$ for all $k$. The $h$-$\mu$ conditions of [23,~Def.~6.4.1], or more precisely 
the $h'$-$\mu$ conditions with $h' = -1$, require only $g \in C^{1,\mu}$ which is 
satisfied by hypothesis. The hypotheses [23, (6.4.2)-(6.4.3)] are also satisfied. 

  Essentially the same manipulations are performed on the boundary system. Consider 
for instance the (most complicated) $3^{\rm rd}$ order boundary operator $\{4\}$. 
One commutes $q^{0\beta}N$ with the Laplacian to obtain schematically an operator 
of the form $\partial(g\partial(g\partial u)) = \partial^{2}(g^{2}\partial u) - 
\partial(g\partial g \partial u)$. The first, leading order, term $g^{2}\partial u$ 
forms one of the boundary operators $B_{rk\gamma}$ in [23,(6.4.15)] with $|\gamma| 
= 2$. The $h'$-$\mu$ conditions on this boundary operator again require 
only $g \in C^{1,\mu}$, which is satisfied. The second term $g\partial g 
\partial u$ forms one of the $g_{r\gamma}$ terms in [23,(6.4.15)] with 
$|\gamma| = 1$. Setting $u = g$, this term is in $L^{1,\hat q}$, for some 
$\hat q > 0$, as required by [23, Thm.~6.4.8]. For this part of $B_{rk\gamma}$, 
one has $h_{r} = 1$ and $p_{r} = 2$, so that $h_{r} + p_{r} \geq 3$, as required 
by [23, (6.4.17)].

  Carrying out the same procedure as needed for the remaining boundary operators 
gives a system of boundary operators with $p_{r} = 0$ for the operators $\{0\}$ 
and $\{1\}$, while $p_{r} = 1, 2$ for the operators $\{3\}$ and $\{4\}$ respectively. 
The terms $h_{r}$ are already defined as following (2.17), and one thus has 
$h_{r} + p_{r} \geq 3$ for all boundary operators, as required by [23, (6.4.17)]. 
As above, it is easily seen that the boundary coefficients satisfy the 
$h'$-$\mu$ conditions. This shows that the hypotheses of [23, Thm.~6.4.8] are 
satisfied, and one concludes that $g \in L^{3,p}$, (assuming 
corresponding regularity in $\widetilde s$ and the boundary data). Given this 
regularity boost, one then iterates this process as needed to obtain $g \in L^{4,p}$ 
or $g \in C^{4,\alpha}$. 

{\endproof}

 Corollary 2.2 leads to the following boundary regularity result.
\begin{theorem} \label{t 2.3.}
  Let $g$ be an AH Einstein metric on a 4-manifold $M$, which admits an 
$L^{2,p}$ conformal compactification $\widetilde g = \rho^{2}g$, $p > 4$, 
with respect to a given background $C^{\infty}$ atlas $\{y^{\mu}\}$ for $M$ 
near $\partial M$, where $\rho = \rho(y^{\mu})$ is an $L^{3,p}(y^{\mu})$ 
defining function. 
 
  If, for a given $m \geq 2$ and $\alpha \in (0,1)$, or $m = \infty$, 
the boundary metric $\gamma  = \widetilde g|_{\partial M}$ is in 
$C^{m,\alpha}(y^{\mu})$, then $g$ admits a $C^{m,\alpha}$ conformal 
compactification $\hat g = \hat \rho^{2}g$, with respect to a 
$C^{m+1,\alpha}$ atlas $\{x^{\mu}\}$ consisting of local boundary 
$\hat g$-harmonic coordinates, with the same boundary metric. Further, 
$\hat \rho = \hat \rho(x^{\mu}) \in C^{m+1, \alpha}(x^{\mu})$. If $\gamma 
\in C^{\omega}(y^{\mu})$, then $\hat g \in C^{\omega}(x^{\mu})$. 

   Moreover, the $x$-coordinates are at least $L^{3,p}(y)$ functions of the 
$y$-coordinates. 
\end{theorem}
{\bf Proof:} 
   Let $\hat g$ be a constant scalar curvature metric conformal to 
$\widetilde g$ on $M$ with $\hat g|_{\partial M} = \gamma$. Thus, for 
$\hat g = u^{2}\widetilde g$, the function $u > 0$ is a solution of the 
Dirichlet problem for the Yamabe equation
\begin{equation}\label{e2.22}
u^{3}\mu  = -6\widetilde \Delta u + \widetilde s u 
\end{equation}
on $M$, with $u = 1$ on $\partial M$ and $\hat s = \mu  = const$. 
It is simplest to choose $\mu  = -1$. Standard methods in elliptic PDE 
then give an $L^{2,p}(y)$ solution to this Dirichlet problem, cf.~[21]. 
Thus, the metric $\hat g$ is $L^{2,p}(y)$ conformally compact, with constant 
scalar curvature. Let $\{x\} = \{x^{\mu}\}$ be a system of local boundary 
$\hat g$-harmonic coordinates near $\partial M$. Then $\hat g \in L^{2,p}(x)$, 
(since harmonic coordinates have optimal regularity), and 
\begin{equation}\label{e2.23}
x \in L^{3,p}(y). 
\end{equation}
Moreover, when restricted to $\partial M$, $x \in C^{m+1,\alpha}(y)$. 
It follows from Corollary 2.2 that $\hat g$ then has the same regularity as the 
boundary metric $\gamma$ in the $x$-coordinates, i.e.~$\hat g \in C^{m,\alpha}(x)$.

 To prove that $\hat \rho \in C^{m+1,\alpha}(x)$, standard formulas relating 
the Ricci curvature of $\hat g$ with that of $g$, cf.~[2, (1.4-(1.5)] or 
[7, Ch.1J], give
$$\hat z \equiv \hat Ric - \frac{\hat s}{4}\hat g = -2\hat \rho^{-1}
(\hat D^{2} \hat \rho - \frac{\hat \Delta \hat \rho}{4}\hat g),$$
so that $\hat D^{2} \hat \rho - \frac{\hat \Delta \hat \rho}{4}\hat g = 
-\frac{1}{2}\hat \rho \hat z$. Now apply the divergence operator $\hat 
\delta$ to both sides of this equation. On the one hand, a simple 
computation gives, (dropping the hats from the notation), $\delta D^{2}f 
= -d\Delta f - Ric(\nabla f)$ and $\delta(f g) = -df$. On the other hand, 
for metrics of constant scalar curvature, $\delta (fz) = f \delta z 
- z(\nabla f) = - z(\nabla f)$, where the last equation follows from 
the contracted Bianchi identity. These calculations then give
$\frac{3}{4}d\hat \Delta \hat \rho = -\frac{\hat s}{4}d\hat \rho 
- \frac{3}{2}\hat z(d\hat \rho)$,
or equivalently
\begin{equation}\label{e2.24}
\hat \Delta d \hat \rho = -\frac{\hat s}{3}d \hat \rho 
- 2\hat z(d \hat \rho).
\end{equation}
Since $\hat s$ is constant and $d$ commutes with $\hat \Delta$, it follows that 
$\hat z(d\hat \rho)$ is exact, so that
\begin{equation}\label{e2.25}
- 2\hat z(d \hat \rho) = d\phi,
\end{equation}
for some function $\phi$. Thus (2.24) is equivalent to 
\begin{equation}\label{e2.26}
\hat \Delta \hat \rho +\frac{\hat s}{3}\hat \rho = \phi,
\end{equation}
(where an undetermined constant has been absorbed into $\phi$). This is an elliptic 
equation for $\hat \rho$, with $\hat \rho = 0$ on $\partial M$, and so one may use 
elliptic boundary regularity results to determine the smoothness of $\hat \rho$. 
To do this, recall that $\hat g \in C^{m,\alpha}(x)$ and $\hat z \in C^{m-2,\alpha}(x)$. 
Suppose first that
\begin{equation}\label{e2.27}
\hat \rho \in C^{k,\alpha}(x),
\end{equation}
for some $k$, $1 \leq k \leq m$. Then $d\hat \rho \in C^{k-1,\alpha}(x)$ and since 
$\hat z \in C^{m-2,\alpha}(x)$, it follows from (2.25) that $d\phi \in C^{\ell,\alpha}$, 
where $\ell = \min(m-2, k-1)$. Hence, $\phi \in C^{\ell+1,\alpha}(x)$. In the $x$-coordinates, 
the Laplacian $\hat \Delta$ has the form $\hat \Delta = \hat g^{\mu\nu}\partial_{x_{\mu}}
\partial_{x_{\nu}}$, and the Schauder elliptic boundary estimates, (cf.~[14] for instance), 
for the equation (2.26) then give 
\begin{equation}\label{e2.28}
\hat \rho \in C^{\ell+3,\alpha}(x),
\end{equation}
provided $\ell+1 \leq m$. This gives an increase in the regularity of $\hat \rho$ by 2 
derivatives over (2.27), and hence by induction it follows that 
\begin{equation}\label{e2.29}
\hat \rho \in C^{m+1,\alpha}(x)
\end{equation}
provided, (for instance), $\hat \rho \in C^{1,\alpha}(x)$. 

  To prove this last statement, note that $\hat \rho = u\rho$, (since $g = 
\hat \rho^{-2}\hat g = \hat \rho^{-2}u^{2}\widetilde g = \hat \rho^{-2}u^{2}\rho^{2}g$). 
One has $\rho \in L^{3,p}(y)$ by assumption and $u\in L^{2,p}(y)$, so that 
$\hat \rho \in L^{2,p}(y)$. Then (2.23) gives $\hat \rho \in L^{2,p}(x) \subset 
C^{1,\alpha'}(x)$, for some $\alpha' > 0$. (The fact that $\alpha'$ may be less than 
$\alpha$ is of no consequence). 

{\endproof}

  One expects that the regularity conclusions in Theorem 2.3 are optimal. Namely, 
it seems unlikely that the regularity of $\widetilde g$ itself can be improved 
without further hypotheses, for example on the scalar curvature $\widetilde s$ 
or on the conformal factor $u$ relating $\hat \rho$ with $\rho$.

\begin{remark} \label{r 2.4.}
  {\rm {\bf (i).} Proposition 2.1, Corollary 2.2 and Theorem 2.3 have all been 
phrased globally. However, the proofs of these results are completely local, 
and so local versions of these results hold equally well. 

  {\bf (ii).} We point out here that the proof of [2,~Thm.2.4] contains 
a small gap. Namely, [2,~Lemma 1.3] does not hold for the mixed components 
$Ric_{0i}$ of the Ricci curvature, when $A \neq  0$. The mixed 
components $Ric_{0i}$ are not determined by $\widetilde s$ and the 
boundary metric, modulo lower order terms as in (2.8) but instead 
are given by (2.15), which is $2^{\rm nd}$ order in the ambient metric. 
Since [2,~Thm.2.4] uses the Yamabe gauge for which $A \neq 0$, one 
does not directly obtain a regularity estimate for $Ric_{0i}$ in this 
gauge. My thanks to Robin Graham and Dylan Helliwell for pointing out 
this gap. 

  The proof of Theorem 2.3 above fixes this gap, via the boundary condition 
$\{4\}$ above. Alternately, it is straightforward to verify that one can also 
prove [2, Thm.2.4] by the same methods used there by adding the boundary 
conditions $\{4\}$. Very briefly, in place of the single Neumann-type boundary 
condition $\{2\}$ used in [2,~Thm.2.4], one uses the pair of Neumann-type 
boundary conditions $\{2\}$ and $\{4\}$, to obtain regularity in the normal 
and mixed directions. The proof of regularity in the tangential directions 
remains the same. 

{\bf (iii).} A version of Theorem 2.3 has been proved in all even dimensions 
recently by Dylan Helliwell, [16]. The proof uses the ideas of the proof 
above in dimension 4, together with the Fefferman-Graham ambient obstruction 
tensor in higher dimensions, in place of the Bach tensor. 
}
\end{remark}

  From certain perspectives, the best compactifications are {\it geodesic} 
compactifications, defined by the property that 
\begin{equation} \label{e2.30}
\bar g = t^{2}g, 
\end{equation}
where $t(x) = dist_{\bar g}(x, \partial M)$. The integral curves of 
$\bar \nabla t$ are then geodesics, orthogonal to $\partial M$ and so the 
Gauss Lemma gives the splitting
\begin{equation} \label{e2.31}
\bar g = dt^{2} + g_{t}, 
\end{equation}
near $\partial M$, where $g_{t}$ may be identified as a curve of 
metrics on $\partial M$ with $g_{0} = \gamma$. Similarly, the metric 
$g$ splits as $g = d\log t^{2} + t^{-2}g_{t}$, so that $r = -\log t$ is a 
geodesic parameter on $(M, g)$. It is well-known that $C^{2}$ conformally 
compact Einstein metrics admit a geodesic compactification, cf.~[12] or [15]. 
Theorem 2.3 gives the following result on the smoothness of the geodesic 
compactification.

\begin{corollary} \label{c 2.5.}
  If $g$ is an $L^{2,p}$ conformally compact Einstein metric on a 
4-manifold M, with $C^{m,\alpha}$ boundary metric $\gamma$, then the 
geodesic compactification $\bar g = t^{2}g$ is $C^{m-1,\alpha}$ smooth, 
in harmonic coordinates. The same result holds with respect to $C^{\infty}$ 
and $C^{\omega}$. 
\end{corollary}
{\bf Proof:} 
By Theorem 2.3, there exists a $C^{m,\alpha}$ compactification $\widetilde g = 
\rho^{2}g$ of $g$ in harmonic coordinates. Writing $t = \omega\rho$, the defining 
equation for $t$, i.e.~$|\bar \nabla t|_{\bar g}^{2} = 1$, is equivalent to
$$2(\widetilde \nabla \rho)\log \omega  + \rho|\widetilde \nabla \log 
\omega|_{\widetilde g}^{2} = \rho^{-1}(1-|\widetilde \nabla \rho|_{\widetilde 
g}^{2}). $$
This is a first order, non-characteristic PDE, with coefficients in 
$C^{m,\alpha}$ and right hand side in $C^{m-1,\alpha}.$ Hence, the 
solution $\omega $ is in $C^{m-1,\alpha}(\bar M)$. 

{\endproof}

 We are now in position to prove Theorem 1.1.

{\bf Proof of Theorem 1.1.}

 We first set up the local Cauchy problem for the Bach equation (2.1). 
As local coordinates, choose geodesic coordinates $(t, x^{i})$ where, 
given a compact metric $g$ on $\bar M$, $t(x) = dist_{g}(x, \partial M)$ 
and $x^{i}$ are local coordinates on $\partial M$ extended into $M$ to be 
invariant under the flow of $T = \nabla t$. Thus, the metric splits in 
these coordinates as in (2.31). (The bar has been dropped from the notation). 
In particular, $g_{0i} = 0$ and $g_{00} = 1$ in these coordinates. Note 
however that the Bach equation (2.1) is {\it not} an elliptic system in 
these coordinates.

 Since (2.1) is a $4^{\rm th}$ order equation, Cauchy data consist of 
prescribing $g$, or equivalently $g_{t}$ in (2.31), and its first three Lie 
derivatives with respect to $T$ at $t = 0$. This data may be freely chosen 
at $\partial M$, but we choose data agreeing with that of the 
Fefferman-Graham expansion \eqref{e1.4} of a conformally compact 
Einstein metric. Thus, set
\begin{equation} \label{e2.32}
g_{(0)} = g_{0} = \gamma , \ \ g_{(1)} = {\mathcal L}_{T}g|_{t=0} = 0, 
\end{equation}
where $\gamma$ is an arbitrary $C^{\omega}$ Riemannian metric on 
$\partial M$. For a conformally compact Einstein metric, the term 
$g_{(2)}$ is intrinsically determined by $\gamma$, (via the Einstein 
equations \eqref{e1.2}), as
\begin{equation} \label{e2.33}
({\mathcal L}_{T}g)^{2}|_{t=0} = g_{(2)} = -\tfrac{1}{2}(Ric_{\gamma} - 
\frac{s_{\gamma}}{4})\gamma . 
\end{equation}
Finally let 
\begin{equation} \label{e2.34}
g_{(3)} = ({\mathcal L}_{T}g)^{3}|_{t=0} = \sigma  
\end{equation}
be an arbitrary transverse-traceless $C^{\omega}$ symmetric bilinear 
form on $(\partial M, \gamma)$, cf.~ again the discussion following 
\eqref{e1.4}. This set of Cauchy data is clearly non-characteristic on 
$\partial M$. We recall that all higher order terms in the expansion 
\eqref{e1.4} are determined by $g_{(0)}$ and $g_{(3)}$. In fact, if one 
defines $g^{k}$ by  $g^{k} = t^{-2}\bar g^{k}$ and 
$$\bar g^{k} \equiv dt^{2} + g_{(0)} + tg_{(1)} + t^{2}g_{(2)} + 
t^{3}g_{(3)} + \cdots + t^{k}g_{(k)},$$
so that $\bar g^{k}$ is a truncation of the Taylor series of $\bar g$, 
then the coefficients $g_{(j)}$ are uniquely determined by the property that
\begin{equation}\label{e2.35}
||Ric_{g^{k}} + 3g^{k}||_{\bar g} = O(t^{k-2}).
\end{equation}
With the exception of $g_{(0)}$ and $g_{(3)}$, one finds that $g_{(j)}$ 
depends on the lower order terms $g_{(l)}$, $l < j$, and their $x$-derivatives 
up to second order, cf.~[12]. 

 Now the system (2.1) has real-analytic coefficients, and the Cauchy 
data above are real-analytic. Of course the boundary $\{t = 0\}$ at 
$\partial M$ is real-analytic in the given coordinates $(t, x^{i})$. 
Hence, the Cauchy-Kovalewsky theorem, cf.~[17], implies there is a 
unique $C^{\omega}$ metric $g$, given in the form (2.31) and defined 
on a thickening $M = [0,\varepsilon)\times \partial M$ of $\partial M$, 
which satisfies the Bach equation (2.1), and satisfies the prescribed 
Cauchy data (2.32)-(2.34).

 Since the curve of metrics $g_{t}$ on $\partial M$ as in \eqref{e2.31} 
is real-analytic in $t$, it is given by its Taylor expansion at $t = 0$. 
Now recall that conformally Einstein metrics are Bach-flat, and so are solutions 
of the equations \eqref{e2.1}. Via the Bach equations, the higher order 
coefficients $g_{(n)}$, $n \geq 4$, in the Taylor expansion of the solution 
$g$ are determined inductively by the lower order terms $g_{(j)}$, 
$0 \leq j \leq 3$ and their $x$-derivatives. Since, by construction in 
(2.32)-(2.34), these lower order terms are determined by the Einstein 
equations, it follows immediately by uniqueness of analytic solutions 
that the higher order terms $g_{(n)}$ are also determined by the Einstein 
equations. Hence the Taylor series of $g_{t}$ is the same as the Fefferman-Graham 
series \eqref{e1.4}. Equivalently, via \eqref{e2.35}, one sees that the 
compactified metric $\bar g$ is conformally Einstein, to infinite order 
at $\partial M$. Analyticity then implies that $\bar g$ is exactly 
conformally Einstein, and moreover that $g = t^{-2}\bar g$ is an AH 
Einstein metric defined near $\partial M$.

 If $g'$ is any other AH Einstein metric with $L^{2,p}$ conformal 
compactification, and with given boundary data $(\gamma, \sigma)$, 
then by Corollary 2.5, the geodesic compactification of $g'$ is 
real-analytic. Hence $g' = g$ up to isometry, so that $g$ is unique among 
AH Einstein metrics (with a weak compactification).
{\endproof}

 As described in [4], the solution to the Einstein equations given by 
Theorem 1.1 can be analytically continued past $N = \partial M$ onto the 
``other side'', to obtain a deSitter-type vacuum solution of the Einstein 
equations. This is a Lorentz metric ${\mathfrak g}$, satisfying the Einstein 
equations with positive cosmological constant, i.e.
\begin{equation} \label{e2.36}
Ric_{\mathfrak g} = 3{\mathfrak g}. 
\end{equation}
This Lorentz metric is $C^{\omega}$ conformally compact, and defined at 
least in the region $M = \partial M\times [0,\varepsilon)$, for some 
$\varepsilon > 0$. Hence, the solution is geodesically complete to the 
future of some Cauchy surface, with real-analytic ${\mathcal I}^{+}$.

 Thus, the Lorentzian version of Theorem 1.1 is the following:

\begin{theorem} \label{t 2.6.}
  Let $N$ be a closed 3-manifold, and let $(\gamma, \sigma)$ be a 
pair consisting of a real-analytic Riemannian metric $\gamma $ on $N$, 
and a real-analytic symmetric bilinear form $\sigma$ on $N$ satisfying 
$\delta_{\gamma}\sigma  = tr_{\gamma}\sigma  = 0$. Then there exists a 
unique vacuum solution to the Einstein equations (2.36) with 
cosmological constant $\Lambda  = 3$, which is $C^{\omega}$ conformally 
compact, defined in a neighborhood of ${\mathcal I}^{+}$, and for which the 
geodesic compactification $\bar {\mathfrak g} = t^{2}{\mathfrak g}$ satisfies
\begin{equation} \label{e2.37}
\bar {\mathfrak g} = (-dt^{2} + \gamma  - t^{2}g_{(2)} - t^{3}\sigma  
+ t^{4}g_{(4)} + ... ). 
\end{equation}
\end{theorem}
{\bf Proof:}
 Given the analyticity from Theorem 1.1 and Corollary 2.5, this is 
proved in [4]. The Fefferman-Graham expansion (1.4) and its basic properties 
holds equally well for Lorentzian deSitter-type vacuum solutions of the 
Einstein equations, cf.~[12]. The terms $g_{(j)}$ in (2.37) are the same as 
those given for the Riemannian AH Einstein metrics in (2.35). Note then that 
formally, the expansion (2.37) is obtained from the expansion (1.5) by replacing 
$t$ by $it$, and dropping any $i$ factors, giving a form of ''Wick rotation" in 
this situation. This is explained in more detail in [4].

  Alternately, one can prove Theorem 2.6 directly, since a Lorentzian vacuum 
solution (2.36) is also Bach-flat. The proof of Theorem 1.1 given above in the 
Riemannian AH setting then carries through in the Lorentzian deSitter-type 
setting in exactly the same way.
{\endproof}

\begin{remark} \label{r 2.7.}
 {\rm  This result gives a simple proof of a result of H. Friedrich [13], 
obtained by solving the conformal Einstein equations, in the special case of 
analytic initial data. A third proof of this result has recently been given 
by A. Rendall [25], using degenerate Fuchsian systems, analogous to the 
original arguments of Fefferman-Graham [12]. }
\end{remark}

\section{Infinitesimal Einstein Deformations and Diffeomorphisms.}
\setcounter{equation}{0}

  This section is a bridge between the previous and next sections. We begin with 
a brief discussion of the Fefferman-Graham expansion [12] in all dimensions and 
then discuss a weak nondegeneracy result from [6] which will be needed in the 
proof of Theorem 1.2. 

  Let $g$ be a conformally compact Einstein metric on a compact $(n+1)$-manifold $M$ 
with boundary $\partial M$ which has a $C^{2}$ geodesic compactification as in 
(2.30). The metric $\bar g$ then splits in geodesic boundary coordinates, as in (2.31):
\begin{equation} \label{e3.1}
\bar g = dt^{2} + g_{t}, 
\end{equation}
near $\partial M$. Each choice of boundary metric $g_{0} = \gamma \in [\gamma]$ 
determines a unique geodesic defining function $t$. Now suppose for the moment that 
the boundary metric $\gamma$ is $C^{\infty}$ smooth. Then by Corollary 2.5 when 
$n = 3$, or by [10] for general $n$, $\bar g$ is $C^{\infty}$ smooth when $n$ is 
odd, and is $C^{\infty}$ polyhomogeneous when $n$ is even. Hence, the curve $g_{t}$ 
has a Taylor-type series in $t$ - the Fefferman-Graham expansion [12]. The exact 
form of the expansion depends on whether $n$ is odd or even. If $n$ is odd, 
one has a power series expansion 
\begin{equation} \label{e3.2}
g_{t} \sim g_{(0)} + t^{2}g_{(2)} + \cdots + t^{n-1}g_{(n-1)} + t^{n}g_{(n)} 
+ \cdots , 
\end{equation}
while if $n$ is even, the series is polyhomogeneous,
\begin{equation} \label{e3.3}
g_{t} \sim g_{(0)} + t^{2}g_{(2)} + \cdots + t^{n}g_{(n)} + t^{n}\log t \ {\mathcal H} + 
\cdots . 
\end{equation}
In both cases, this expansion is even in powers of $t$, up to $t^{n}$. The coefficients 
$g_{(2k)}$, $k \leq [n/2]$, as well as the coefficient ${\mathcal H}$ when $n$ is even, 
are explicitly determined by the boundary metric $\gamma  = g_{(0)}$ and the Einstein 
condition (1.2), cf.~[11], [12]. For $n$ even, the series (3.3) has terms of the form 
$t^{n+k}(\log t)^{m}$. 

  For any $n$, the divergence and trace (with respect to $g_{(0)} = \gamma$) of $g_{(n)}$ 
are determined by the boundary metric $\gamma$; in fact there is a symmetric bilinear form 
$r_{(n)}$ and scalar function $a_{(n)}$, both depending only on $\gamma$ and its derivatives 
up to order $n$, such that 
\begin{equation} \label{e3.4}
\delta_{\gamma}(g_{(n)} + r_{(n)}) = 0, \ \ {\rm and} \ \ tr_{\gamma}(g_{(n)} + r_{(n)}) = 
a_{(n)}.
\end{equation}
For $n$ odd, $r_{(n)} = a_{(n)} = 0$. However, beyond the relations (3.4), the term 
$g_{(n)}$ is not determined by $g_{(0)}$; it depends on the ``global'' structure of the 
metric $g$. The higher order coefficients $g_{(k)}$ of $t^{k}$ and coefficients $h_{(km)}$ 
of $t^{n+k}(\log t)^{m}$, are then determined by $g_{(0)}$ and $g_{(n)}$ via the Einstein 
equations. The equations (3.4) are constraint equations, and arise from the Gauss-Codazzi 
and Gauss and Riccati equations on the level sets $S(t) = \{x: t(x) = t\}$ in the limit 
$t \rightarrow 0$. 

   Now suppose $k$ is an infinitesimal Einstein deformation of $(M, g)$, so that 
$k$ satisfies
\begin{equation} \label{e3.5}
L_{E}(k) \equiv 2\frac{d}{ds}(Ric_{g+sk} + n(g+sk)) = D^{*}Dk - 2R(k) - 
2\delta^{*}\beta (k) = 0,
\end{equation}
where $\beta$ is the Bianchi operator $\beta(k) = \delta k + \frac{1}{2}dtr k$. 
Suppose for the moment that $k$ is $C^{\infty}$ polyhomogeneous smooth up to 
$\partial M$ and preserves the geodesic boundary coordinates near $\partial M$, 
so that $k_{0\alpha} = 0$. If 
\begin{equation}\label{e3.6}
k = O(t),
\end{equation}
on approach to $\partial M$, then the discussion above on the Fefferman-Graham 
expansion implies the stronger decay
\begin{equation}\label{e3.7}
k = O(t^{n}).
\end{equation}
If moreover one assumes the stronger condition that
\begin{equation}\label{e3.8}
k = o(t^{n}),
\end{equation}
then the induced variation of the terms $g_{(0)}$ and $g_{(n)}$ in 
\eqref{e3.2}-\eqref{e3.3} vanishes and, again in view of the discussion on the 
expansions above, one has 
\begin{equation}\label{e3.9}
k = o(t^{\nu}),
\end{equation}
for all $\nu < \infty$. In this situation, one would expect that $k \equiv 0$ 
near $\partial M$. More generally, if $k$ is as above but is not necessarily in 
geodesic gauge, then near $\partial M$, $k$ should be a ``pure gauge'' deformation, 
i.e.~$k = \delta^{*}Z$, for some vector field $Z$ on $M$ with $X = 0$ on $\partial M$. 
These expectations do in fact hold, and are proved in [6]; this corresponds 
to a unique continuation property at infinity for AH Einstein metrics and 
their linearizations. 

   For the work to follow in \S 4, we need to discuss this in somewhat more 
detail. Thus, suppose $(M, g)$ is a $C^{2,\alpha}$ conformally compact Einstein 
metric. In view of \eqref{e3.5}, the simplest gauge choice for infinitesimal 
Einstein deformations $\kappa$ of $g$ is the Bianchi gauge
\begin{equation}\label{e3.10}
\beta(\kappa) = 0.
\end{equation}
In this case, $\kappa$ satisfies the elliptic equation
\begin{equation}\label{e3.11}
L(\kappa) = D^{*}D\kappa - 2R(\kappa) = 0.
\end{equation}
A well-known result of Biquard [8] also gives a converse to this statement. 
Namely, if $\kappa$ is a solution of \eqref{e3.11} satisfying \eqref{e3.6}, 
then \eqref{e3.10} holds. In fact, one then has
\begin{equation}\label{e3.12}
\delta \kappa = 0 \ \ {\rm and} \ \ tr \kappa = 0.
\end{equation}
To prove \eqref{e3.12}, the trace of \eqref{e3.11} gives the equation
$$-\Delta tr \kappa + 2n tr \kappa = 0.$$
Since $|tr \kappa| \rightarrow 0$ at infinity, it follows immediately from 
the maximum principle that $tr \kappa = 0$. Combining this with \eqref{e3.10} 
shows that $\delta \kappa = 0$ also.  

   We also note the well-known fact, proved via elliptic regularity in 
weighted H\"older spaces associated to the equation \eqref{e3.11}, that for 
$(M, g)$ as above, if $\kappa$ satisfies \eqref{e3.11} and either $\kappa \in 
L^{2}(M, g)$ or \eqref{e3.6} holds for $\kappa$, then \eqref{e3.7} holds, i.e.
\begin{equation}\label{e3.13}
\kappa = O(t^{n}).
\end{equation}
In addition, an analysis of the behavior of the indicial roots of \eqref{e3.11} 
shows that one also has
\begin{equation}\label{e3.14}
\kappa (N,\cdot) = O(t^{n+1}),
\end{equation}
where $N = \rho \partial_{\rho}$, with $\rho$ the given defining function 
for $\partial M$ in $M$; the decay estimates \eqref{e3.13}-\eqref{e3.14} are 
proved in [8], [20] or [22].

  Given this background, the following result is proved in [6,~Cor.~4.4], 
and will be used in the proof of Theorem 4.1.

\begin{proposition} \label{p3.1}
Let $g$ be a $C^{3,\alpha}$ conformally compact Einstein metric on $M$, and 
suppose
\begin{equation}\label{e3.15}
\pi_{1}(M, \partial M) = 0.
\end{equation}
If $\kappa$ is a symmetric bilinear form in $L^{2}(M, g)$ satisfying 
\eqref{e3.11} and \eqref{e3.8} holds in a neighborhood of $\partial M$, 
then 
\begin{equation} \label{e3.16}
\kappa \equiv 0 ,
\end{equation}
on $(M, g)$. 
\end{proposition}

\begin{remark} \label{r3.2} 
{\rm  Proposition 3.1 proves a weak nondegeneracy property conjectured in [22]. 
As noted above, it corresponds to a unique continuation property at infinity 
for solutions of the linearized AH Einstein equations. A local version of this 
result also holds, (where \eqref{e3.8} holds only on approach to a portion of 
the boundary); this will not be used here however. The topological condition 
\eqref{e3.15} is needed to ensure that the (iterative) use of the local unique 
continuation property for solutions of \eqref{e3.11} extends consistently to 
give a global uniqueness on the full manifold $M$. 

  We expect that Proposition 3.1 is false in general if the assumption 
$\pi_{1}(M, \partial M) = 0$ is dropped, for example if $\partial M$ is 
not connected. However, this is not known and it would be of interest to 
find some concrete counterexamples. 

  On the other hand, if $\partial M$ is connected and $(M, g)$ has no local Killing 
fields, (i.e.~there are no Killing fields on the universal cover $\widetilde M$ of 
$M$), then the proof Proposition 3.1 in [6] holds without the assumption 
\eqref{e3.15}. In particular, it follows that Proposition 3.1 holds for generic 
AH Einstein metrics on $M$ provided $\partial M$ is connected.

   Finally for the application in \S 4, we note that the pointwise assumption 
\eqref{e3.8} in Proposition 3.1 may be weakened to the analogous assumption on 
the $L^{2}$ norm of $\kappa$ over the spheres $S(t)$ as $t \rightarrow 0$. This 
follows again from elliptic regularity associated with the equation \eqref{e3.11}, 
cf.~[6, Lemma 4.2] for further details. }
\end{remark}

\section{The Banach Manifold ${\mathcal E}_{AH}$.}
\setcounter{equation}{0}

 In this section we prove that the moduli space of AH Einstein metrics 
on a given $(n+1)$-manifold is naturally an infinite dimensional Banach 
manifold, assuming it is non-empty. This is essentially the content of 
Theorem 1.2, but the full version is proved in \S 5. The work 
in this section uses the methods developed by Graham-Lee [15] and 
Biquard [10], as well as the work of White [26], [27].

 We begin by describing the function spaces to be used. First, let 
$\rho_{0}$ be a fixed $C^{\omega}$ defining function for $\partial M$ 
in $M$. Throughout \S 4, the defining function $\rho_{0}$ will be kept 
fixed and only compactifications with respect to $\rho_{0}$, will be considered, 
i.e.
\begin{equation} \label{e4.1}
\widetilde g = \rho_{0}^{2}\cdot  g. 
\end{equation}
The situation where $\rho_{0}$ varies over the family of smooth defining 
functions is discussed in \S 5. Given $\rho_{0}$, define the 
function $r = r(\rho_{0})$ on $M$ by
\begin{equation} \label{e4.2}
r = - \log(\tfrac{\rho_{0}}{2}). 
\end{equation}

 Let $Met^{m,\alpha}(\partial M)$ be the space of $C^{m,\alpha}$ 
Riemannian metrics on $\partial M,$ so that $Met^{m,\alpha}$ is an open 
cone in the Banach space ${\mathbb S}^{m,\alpha}(\partial M)$ of symmetric 
bilinear forms on $\partial M$. The space $Met^{m,\alpha}(\partial M)$ 
is given the $C^{m,\alpha'}$ topology, for a fixed $\alpha' < \alpha$, 
so that bounded sequences in the $C^{m,\alpha}$ norm have convergent 
subsequences. In this topology, $Met^{m,\alpha}$ is separable, cf.~[26]. 
Next let ${\mathbb S}^{k,\beta}(M)$ be the Banach space of $C^{k,\beta}$ 
symmetric bilinear forms on $M$, and let ${\mathbb S}^{k,\beta}(\bar M)$ 
be the corresponding space of forms on the closure $\bar M$, again with 
the $C^{k,\beta'}$ topology, $\beta' < \beta$. 

 Forms in ${\mathbb S}^{k,\beta}(M)$ have no control or restriction on 
their behavior on approach to $\partial M,$ while those in ${\mathbb 
S}^{k,\beta}(\bar M)$ of course by definition extend $C^{k,\beta}$ up 
to $\partial M.$ Thus, $(m,\alpha )$ determines the regularity of the 
boundary data, while $(k,\beta )$ determines the regularity in the 
interior $M$. These are not necessarily related, unless one has boundary 
regularity results, i.e. regularity of the data up to and including the 
boundary. We will always assume that $m+\alpha \geq  k+\beta$, and $k \geq 2$, 
$\alpha, \beta \in (0,1)$. 

 Let $g$ be a complete Riemannian metric of bounded geometry on $M$, i.e.~$g$ 
has bounded sectional curvature and injectivity radius bounded below on $M$. 
Following [15] and [8], define the weighted H\"older spaces 
${\mathbb S}_{\delta}^{k,\beta}(M) = {\mathbb S}_{\delta}^{k,\beta}(M, g)$ to 
be the Banach space of symmetric bilinear forms $h$ on $M$ such that
\begin{equation} \label{e4.3}
h = e^{-\delta r}h_{0}, 
\end{equation}
where $h_{0}\in{\mathbb S}^{k,\beta}(M)$ satisfies 
$||h_{0}||_{C^{k,\beta}(M, g)} \leq C$, for some constant $C <  \infty$. 
Here the norm is the usual $C^{k,\beta}$ norm with respect to the metric $g$, 
given by
\begin{equation} \label{e4.4}
||h_{0}||_{C^{k,\beta}(M, g)} = \sum_{j\leq 
k}||\nabla^{j}h_{0}||_{C^{0}(M, g)} + ||\nabla^{k}h_{0}||_{C^{\beta}(M, g)}. 
\end{equation}
Thus $h\in{\mathbb S}_{\delta}^{k,\beta}(M)$ implies that $h$ and its 
derivatives up to order $k$ with respect to $g$ decay as $e^{-\delta r}$ as 
$r \rightarrow  \infty$. The weighted norm of $h$ is then defined as 
\begin{equation} \label{e4.5}
||h||_{C_{\delta}^{k,\beta}(M)}= ||h_{0}||_{C^{k,\beta}(M)}. 
\end{equation}
The norms in (4.4) and (4.5) depend only on $C^{k,\beta}$ 
quasi-isometry class of $g$; two metrics $g$ and $g'$ are $C^{k,\beta}$ 
quasi-isometric if, in a fixed local coordinate system, the linear map 
$g^{-1}g' $ is bounded away from 0 and $\infty $ in $C^{k,\beta}(M).$ 
Hence the spaces ${\mathbb S}_{\delta}^{k,\beta}(M)$ depend only on the 
$C^{k,\beta}$ quasi-isometry class of $g$.

 Now suppose the metric $g$ is conformally compact, with 
compactification $\widetilde g$ as in (4.1). One may then define the 
$C^{k,\beta}$ norm of $h$ above also with respect to $\widetilde g$. Using 
standard formulas for conformal changes of metric gives, for any $j$, 
$\beta \geq  0$,
\begin{equation} \label{e4.6}
||\widetilde \nabla^{j}h||_{C^{\beta}(\widetilde g)} = 
||\rho_{0}^{-2-j-\beta}\nabla^{j}h||_{C^{\beta}(g)} + \ \ {\rm lower \ order \ terms}. 
\end{equation}

 Given these preliminaries, one can construct a natural or ``standard'' AH metric 
associated to any boundary metric $\gamma\in Met^{m,\alpha}(\partial M)$. 
This is first done in a collar neighborhood $U$ of $\partial M$ on 
which $d\rho_{0}$ is non-zero, and then later extended to a 
metric on all of $M$. Choose a fixed identification of $U$ with 
$I\times \partial M$ so that $\rho_{0}$ corresponds to the variable on $I$. 
Recalling that $\rho_{0}$ is fixed, define the $C^{m,\alpha}$ 
hyperbolic cone metric $g_{U} = g_{U}(\gamma , \rho_{0})$ in $U$ by
\begin{equation} \label{e4.7}
g_{U} = dr^{2} + \sinh^{2}r\cdot \gamma , 
\end{equation}
for $r$ as in (4.2). Observe that the dependence of $g_{U}$ is 
$C^{\omega}$ in $\gamma$, (and also in $\rho_{0}$). Also if 
$\gamma_{1}$ and $\gamma_{2}$ are $C^{m,\alpha}$ quasi isometric 
boundary metrics, then $g_{U}(\gamma_{1})$ and $g_{U}(\gamma_{2})$ are 
$C^{m,\alpha}$ quasi isometric.

 If $\{e_{i}\}$ is a local orthonormal frame for $g_{U},$ with $e_{1} = 
\partial_{r}$, then one easily verifies that the sectional curvatures 
$K_{ij}$ of $g_{U}$ in the direction $(e_{i}, e_{j})$ are given by
$$K_{1i} = - 1, K_{jk} = \frac{1}{\sinh^{2}r}(K_{\gamma})_{jk} -  
\coth^{2}r, $$
where $i,j,k$ run from $2$ to $n+1$ and $K_{\gamma}$ is the sectional 
curvature of $\gamma$. This implies that the curvature of $g_{U}$ 
decays to that of the hyperbolic space $H^{n+1}(- 1)$ at a rate of 
$O(e^{- 2r}) = O(\rho_{0}^{2})$. The same decay holds for the covariant 
derivatives of the curvature, up to order $m- 2+\alpha$. In particular 
by (4.3)-(4.5)
\begin{equation} \label{e4.8}
Ric_{g_{U}} + n\cdot  g_{U} \in  {\mathbb S}_{2}^{m- 2,\alpha}(M). 
\end{equation}
The metric $g_{U}$ is $C^{m,\alpha}$ conformally compact. In fact if 
$\widetilde g_{U}$ is the compactification (4.1) of $g_{U}$, then a simple 
computation gives
$$\widetilde g_{U} = d\rho_{0}^{2} + (1-\tfrac{1}{4}\rho_{0}^{2})^{2}\cdot 
\gamma .$$ 

 The metric $g_{U}$ will be viewed as a background metric with which to 
compare other conformally compact metrics with the same boundary 
metric. Thus suppose $g'$ is any conformally compact metric on $M$, with 
compactification $\widetilde g'$ as in (4.1). Then one may write
\begin{equation} \label{e4.9}
g'|_{U} = g_{U} + h, 
\end{equation}
and we will assume that $h \in  {\mathbb S}^{k,\beta}(M)$. This implies 
that
\begin{equation} \label{e4.10}
\widetilde g'|_{U} = d\rho_{0}^{2} + (1-\tfrac{1}{4}\rho_{0}^{2})^{2}\cdot 
\gamma  + \rho_{0}^{2}h, 
\end{equation}
so that if $\rho_{0}^{2}|h| \rightarrow 0$ on approach to $\partial M$, 
then $\widetilde g'$ is a $C^{0}$ compactification of $g'|_{U}$ with 
boundary metric $\gamma$; here $|h|$ is the pointwise norm of $h$ 
with respect to any smooth metric on $\bar M$. The compactification 
$\widetilde g'$  is $C^{k,\beta}$ when $\rho_{0}^{2}h\in{\mathbb S}^{k,\beta}(\bar M)$. 
Using the relations (4.3)-(4.6), observe that
\begin{equation} \label{e4.11}
h\in{\mathbb S}_{\delta}^{k,\beta}(M) \ \ {\rm with} \ \ \delta  \geq  k+\beta  
\Rightarrow \widetilde h = \rho_{0}^{2}h\in{\mathbb S}^{k,\beta}(\bar M). 
\end{equation}
However, if $\delta < k+\beta$ and $h\in{\mathbb S}_{\delta}^{k,\beta}(M)$, 
then in general, i.e.~without further restrictions, $\rho_{0}^{2}h$ will 
not be in ${\mathbb S}^{k,\beta}(\bar M)$; this is essentially the issue of 
boundary regularity, and will be discussed at the end of \S 4 and in \S 5.

 The standard metrics $g_{U}$ may be naturally extended to all of $M$ 
as follows. Let $\eta = \eta(r)$ be a fixed cutoff function on $M$, with 
$\eta  \equiv 1$ on $U$, $\eta  \equiv 0$ on $M\setminus U'$, where $U'$ 
is a thickening of $U$ on which $d\eta$ is also non-vanishing. If $g_{C}$ 
is any smooth Riemannian metric on the compact manifold $M\setminus U$, 
(so $g_{C}$ is incomplete), then define
\begin{equation} \label{e4.12}
g_{\gamma} = \eta g_{U} + (1 - \eta) g_{C}. 
\end{equation}
Thus for any $\gamma\in Met^{m,\alpha}(\partial M)$, $g_{\gamma}$ in 
(4.12) gives a standard AH metric on $M$, with boundary metric $\gamma$. 
The metric $g_{\gamma}$ on $M$ again depends smoothly on $\gamma $ and 
the choices of the compact metric $g_{C}$ and cutoff $\eta .$ As with 
$\rho_{0}$, we fix the metric $g_{C}$ and cutoff $\eta$ once for all. 
With this understood, one thus has a $C^{\omega}$ smooth (addition) map
\begin{equation} \label{e4.13}
A: Met^{m,\alpha}(\partial M)\times {\mathcal U}_{\delta}^{k,\beta} \rightarrow  
Met_{\delta}^{k,\beta}(M),  A(\gamma , h) = g \equiv  g_{\gamma}+h, 
\end{equation}
where ${\mathcal U}_{\delta}^{k,\beta}$ is the open subset of ${\mathbb 
S}_{\delta}^{k,\beta}(M),$ consisting of those $h$ such that 
$g_{\gamma}+h$ is a well-defined metric on $M$. 

 In view of the decay rate (4.8), the most natural choice of $\delta$ 
is
\begin{equation} \label{e4.14}
\delta  = 2, 
\end{equation}
and we fix this choice for the remainder of this section. The map $A$ is 
clearly injective and the asymptotically hyperbolic (AH) metrics (of 
weight $\delta  = 2$) are defined to be the image of $A$;
\begin{equation} \label{e4.15}
Met_{AH}^{k,\beta} = Im A. 
\end{equation}
The inverse map to $A$, $S: Met_{AH}^{k,\beta} \rightarrow  
Met^{m,\alpha}(\partial M)\times {\mathcal U}_{2}^{k,\beta}$ gives the splitting 
of the AH metric $g$ into its components $g_{\gamma}$ and $h$. Let 
\begin{equation} \label{e4.16}
E_{AH}^{k,\beta} \subset  Met_{AH}^{k,\beta} 
\end{equation}
be the subset of AH Einstein metrics, with topology induced as a subset 
of the product topology. Note that, as discussed in (4.11), metrics in 
$E_{AH}^{k,\beta}$ are $C^{2}$ conformally compact, but not necessarily 
$C^{k,\beta}$ conformally compact with respect to $\rho_{0}$, for $k+\beta > 
2$. Of course Einstein metrics are $C^{\omega}$ in local harmonic 
coordinates, and so $(k,\beta)$ only serves to denote the ambient space 
$Met_{AH}^{k,\beta}$ in which $E_{AH}$ is embedded. 

 Now let $g_{0}$ be a fixed (but arbitrary) background metric in $Met_{AH}^{k,\beta}$ 
with boundary metric $\gamma_{0}$. For $\gamma \in Met^{m,\alpha}(\partial M)$ 
close to $\gamma_{0}$, let 
\begin{equation} \label{e4.17}
g(\gamma) = g_{0} + \eta(g_{\gamma} - g_{\gamma_{0}}).
\end{equation}
Any metric $g \in Met_{AH}^{k,\beta}$ with boundary metric $\gamma$ thus 
has the form $g = g(\gamma) + h$, for $h \in {\mathbb S}_{2}^{k,\beta}$. 
Essentially as in [8], for any $k \geq 2$, define 
\begin{equation}\label{e4.18}
\Phi = \Phi^{g_{0}}: Met_{AH}^{k,\beta} \rightarrow  {\mathbb S}_{2}^{k-2,\beta}(M),
\end{equation}
$$\Phi (g) = \Phi (g(\gamma) + h) = Ric_{g} + ng + (\delta_{g})^{*}\beta_{g(\gamma)}(g),$$
where $\beta_{g(\gamma)}$ is the Bianchi operator with respect to $g(\gamma)$, 
(cf.~\eqref{e3.5}), 
\begin{equation} \label{e4.19}
\beta_{g(\gamma)}(g) = \delta_{g(\gamma)}g + \tfrac{1}{2}d(tr_{g(\gamma)}g). 
\end{equation}
Observe that $\Phi$ is well-defined, by (4.8), (4.14) and the fact 
that $h = g - g(\gamma) \in {\mathbb S}_{2}^{k,\beta}$. Clearly $\Phi$ is 
$C^{\omega}$ in $g$. 

  There are several natural reasons for considering the operator $\Phi$. First, 
it is proved in [8, Lemma I.1.4] that
\begin{equation} \label{e4.20}
Z_{AH}^{k,\beta} \equiv  \Phi^{-1}(0)\cap\{Ric <  0\} \subset  
E_{AH}^{k,\beta}, 
\end{equation}
where $\{Ric < 0\}$ is the open set of metrics with negative Ricci 
curvature. (Here one uses the fact that $\beta_{g(\gamma)}(g) \in{\mathbb S}_{2}^{k-1,\beta}$). 
Further, if $g$ is an AH Einstein metric, i.e. $Ric_{g} = - ng$, with boundary 
metric $\gamma$, which is close to $g_{0}$ and which satisfies $\Phi (g) = 0$, then
\begin{equation} \label{e4.21}
\beta_{g(\gamma)}(g) = 0. 
\end{equation}
As discussed later, the condition (4.21) defines the tangent space of a 
slice to the action of the diffeomorphism group on $Met_{AH}$ and 
$E_{AH}$. Thus, for any $g\in E_{AH}^{k,\beta}$ near $g_{0}$, there exists a 
diffeomorphism $\phi $ such that $\phi^{*}g\in Z_{AH}^{k,\beta}$, cf. 
(4.38) below. Hence, $E_{AH}^{k,\beta}$ differs from $Z_{AH}^{k,\beta}$ 
just by the action of diffeomorphisms.

 Second, as discussed in \S 3, the linearization of the Einstein operator $Ric_{g} + ng$ 
at an Einstein metric $g$ is given by 
\begin{equation} \label{e4.22}
\tfrac{1}{2}D^{*}D -  R - \delta^{*}\beta, 
\end{equation}
acting on the space of symmetric 2-tensors ${\mathbb S}(M)$ on $M$, cf.~[7]. 
The kernel of the elliptic self-adjoint linear operator
\begin{equation} \label{e4.23}
L = \tfrac{1}{2}D^{*}D -  R 
\end{equation}
corresponds to the space of non-trivial infinitesimal Einstein deformations 
in Bianchi-free gauge, analogous to the Jacobi fields for geodesics. An AH 
Einstein metric $g$ on $M$ is called {\it non-degenerate} if 
\begin{equation} \label{e4.24}
K = L^{2}-Ker L = 0, 
\end{equation}
i.e.~if there are no non-trivial infinitesimal Einstein deformations of 
$g$ in $L^{2}(M, g)$. Einstein metrics are critical points of the 
Einstein-Hilbert functional or action, and this corresponds formally to 
the condition that the critical point be non-degenerate, in the sense 
of Morse theory. Recall from \eqref{e3.12} that elements 
in $K$ are transverse-traceless. 

 Now the linearization of $\Phi$ at $g_{0} \in Met_{AH}^{k,\beta}$ 
with respect to the $2^{\rm nd}$ variable $h$ has the simple form
\begin{equation} \label{e4.25}
(D_{2}\Phi)_{g_{0}}(\dot h) = \tfrac{1}{2}D^{*}D \dot h + 
\tfrac{1}{2}(Ric_{g_{0}}\circ \dot h + \dot h \circ Ric_{g_{0}} + 2n \dot h) 
-  R_{g_{0}}(\dot h); 
\end{equation}
this is due to cancellation of the variation of the term 
$(\delta_{g})^{*}\beta_{g(\gamma)}(g)$ with the variation of the Ricci 
curvature, cf. [8, (1.9)]. Hence, if $g_{0}$ is Einstein, then
\begin{equation} \label{4.26}
(D_{2}\Phi)_{g_{0}} = L = \tfrac{1}{2}D^{*}D -  R.
\end{equation}
The variation of $\Phi$ at $g_{0}$ with respect to the $1^{\rm st}$ variable 
$g(\gamma)$ has the form 
$$(D_{1}\Phi)_{g_{0}}(\dot g(\gamma)) = (D_{2}\Phi)_{g_{0}}(\dot g(\gamma)) 
- \delta_{g_{0}}^{*}\beta_{g_{0}}(\dot g(\gamma)),$$ 
where $(D_{2}\Phi)_{g_{0}}(\dot g(\gamma))$ is given by (4.25) with 
$\dot g(\gamma)$ in place of $\dot h$. 

\medskip

 The main result of this section, which leads to a version of Theorem 
1.2 is the following:
\begin{theorem} \label{t 4.1.}
  Suppose $\pi_{1}(M, \partial M) = 0$. At any metric $g\in E_{AH}^{k,\beta}$ 
which is $C^{3,\alpha}$ conformally compact, the map $\Phi$ is a submersion, 
i.e.~the derivative
\begin{equation} \label{e4.27}
(D\Phi )_{g}: T_{g}Met_{AH}^{k,\beta} \rightarrow  
T_{\Phi(g)}{\mathbb S}_{2}^{k-2,\beta}(M) 
\end{equation}
is surjective and its kernel splits in $T_{g}Met_{AH}^{k,\beta}$. 
\end{theorem}
{\bf Proof:}
 By (4.13) and (4.15), $T_{g}Met_{AH} = T_{\gamma}Met^{m,\alpha}(\partial M)
\oplus T_{h}{\mathbb S}_{2}^{k-2,\beta}(M)$. With respect to this splitting, (4.26) 
shows that the derivative of $\Phi$ with respect to the second (i.e. $h$) factor 
is given by 
$$(D_{2}\Phi)_{g} = \tfrac{1}{2}D^{*}D -  R : {\mathbb S}_{2}^{k,\beta}(M) 
\rightarrow  {\mathbb S}_{2}^{k-2,\beta}(M). $$
Now by [8, Prop. I.3.5], $(D_{2}\Phi)_{g}$ is a Fredholm operator 
whose kernel on ${\mathbb S}_{2}^{k,\beta}(M)$ equals the $L^{2}$ kernel $K$ in 
(4.24). Since $(D_{2}\Phi)_{g}$ is self-adjoint on $L^{2}$, 
it has Fredholm index 0, and the cokernel of $(D_{2}\Phi )_{g}$ is 
naturally identified with $K$ in ${\mathbb S}_{2}^{k-2,\beta}(M)$. Thus to 
prove $D\Phi_{g}$ is surjective, it suffices to show that for any 
non-zero $L^{2}$ infinitesimal Einstein deformation $\kappa\in K$, 
there is a tangent vector $X\in T_{g}Met_{AH}$ such that
\begin{equation} \label{e4.28}
\int_{M}\langle (D\Phi)_{g}(X), \kappa \rangle dV_{g} \neq  0. 
\end{equation}

 To do this, let $X = \dot g(\gamma)$, so that $X$ corresponds to a variation 
of the boundary metric $\gamma$ of $g$. Then $D\Phi_{g_{0}}(X)$ has 
the form 
\begin{equation} \label{e4.29}
(D\Phi_{g})(X) = \tfrac{1}{2}D^{*}D \dot g(\gamma) -  R(\dot g(\gamma)) - 
\delta^{*}\beta(\dot g(\gamma)). 
\end{equation}
Let $\gamma$ be the boundary metric induced by $\widetilde g$ in (4.1) on $\partial M$. 
For the following computation, it is convenient to work with the geodesic defining 
function $t$ determined by $\gamma$. Set $r = -\log \frac{t}{2}$, as in (4.2) and 
let $B(r)$ be the $r$-sublevel set of the function $r$ with $S(r) = \partial B(r)$ 
the $r$-level set. We apply the divergence theorem to the integral (4.28) over 
$B(r)$; twice for the Laplacian term in (4.29) and once for the $\delta^{*}$ term. 
Since $\kappa \in Ker L$ and $\delta \kappa = 0$ by \eqref{e3.12}, it follows that 
the integral (4.28) reduces to an integral over the boundary, and gives
\begin{equation} \label{e4.30}
\int_{B(r)}\langle (D\Phi )_{g}(X), \kappa \rangle dV_{g} = 
\int_{S(r)}(\langle \dot g(\gamma), \nabla_{N}\kappa \rangle  -  
\langle \nabla_{N}\dot g(\gamma), \kappa \rangle  - 
\langle \beta(\dot g(\gamma)), \kappa(N) \rangle )dV_{S(r)}, 
\end{equation}
where $N = \nabla r = -t\partial_{t}$ is the unit outward normal.

 To estimate the boundary integrals, the volume form of $S(r)$ satisfies
$$dV_{S(r)} = t^{-n}dV_{\gamma} + O(t^{2}),$$
where $dV_{\gamma}$ is the volume form of the boundary metric. Let 
$\widetilde \kappa = t^{-n}\kappa$. By \eqref{e3.13}, 
$|\widetilde \kappa|_{g}|_{S(r)}$ is uniformly bounded. Setting 
$\hat \kappa = t^{2}\widetilde \kappa$, one has $|\hat \kappa|_{\bar g} = 
|\widetilde \kappa|_{g}$, and so the same is true for $|\hat \kappa|_{\bar g}$. 
A simple calculation from (4.17) gives 
$$\dot g(\gamma) = \eta (\sinh^{2}r)\dot \gamma, \ \nabla_{N}\dot g(\gamma) 
= O(t),$$
so that $|\dot g(\gamma)|_{g} \sim 1$ and $|\nabla_{N}\dot g(\gamma)|_{g} \sim O(t)$ 
as $t \rightarrow 0$. Hence, 
\begin{equation}\label{e4.31}
(\langle \dot g(\gamma), \nabla_{N}\kappa \rangle_{g}  -  
\langle \nabla_{N}\dot g(\gamma), \kappa \rangle)_{g}dV_{S(r)} = 
t^{2}\langle \nabla_{N}\kappa , \dot \gamma \rangle_{\gamma}dV_{S(r)} + O(t)
\end{equation}
$$= \langle \nabla_{N}\hat \kappa - (n-2)\hat\kappa, \dot \gamma 
\rangle_{\gamma} dV_{\gamma} + O(t).$$
By \eqref{e3.14}, $\kappa (N) = O(t^{n+1})$, and hence the last term in (4.30) 
vanishes in the limit $r \rightarrow \infty$. 

   It follows that if (4.28) vanishes in the limit $r \rightarrow 
\infty$, for {\it all} variations $\dot \gamma$, then one must have
\begin{equation} \label{e4.32}
\nabla_{N} \hat \kappa - (n-2)\hat \kappa = o(1),
\end{equation}
weakly, as forms on $T(S(r))$ with respect to $\bar g$. Here $\nabla$ is the 
covariant derivative with respect to $g$, not $\bar g$. Pairing this with the 
bounded form $\hat \kappa$, and using \eqref{e3.14} again, one easily sees that 
(4.32) implies that $\frac{1}{2}N|\hat \kappa|^{2} - n|\hat \kappa|^{2} = o(1)$, 
where the norms are with respect to $\bar g$. Integrating this with respect to the 
volume form on $(S(t), \bar g)$ and using the fact that $\frac{d}{dt}dV_{S(t)} 
= O(t)$, it follows that
\begin{equation} \label{e4.33}
{\tfrac{1}{2}}N\int_{S(t)}|\hat \kappa|^{2}dV_{\gamma} - 
n\int_{S(t)}|\hat \kappa|^{2}dV_{\gamma} = o(1),
\end{equation}
as $t \rightarrow 0$. Using the fact that $N = -t\partial t$, an elementary 
integration then implies that $\int_{S(t)}|\hat \kappa|^{2}dV_{\gamma} 
\rightarrow 0$, and hence $\widetilde \kappa = o(t^{n})$ weakly. However, 
under the assumption $\pi_{1}(M, \partial M) = 0$, this contradicts 
Proposition 3.1, (cf.~also Remark 3.2), which thus proves that (4.28) holds. 

\medskip

 To prove that the kernel of $D\Phi_{g}$ splits, i.e.~it admits a 
closed complement in $T_{g}Met_{AH},$ it suffices to exhibit a bounded 
linear projection $P$ mapping $T_{g}Met_{AH}$ onto $Ker(D\Phi_{g}).$ We 
do this following [27]. Thus, one has 
\begin{equation} \label{e4.34}
KerD\Phi_{g} = \{(\dot \gamma, \dot h):D_{1}\Phi (\dot \gamma) + 
D_{2}\Phi (\dot h) = 0\}.  
\end{equation}
From (4.26), $D_{2}\Phi  = L$ and Im $L = K^{\perp},$ for $K$ as in 
(4.24). Hence $D_{1}\Phi (\dot \gamma)\in  K^{\perp},$ so that 
$\pi_{K}D_{1}\Phi (\dot \gamma) = 0$, i.e. $\dot \gamma\in 
Ker(\pi_{K}D_{1}\Phi)$, where $\pi_{K}$ is orthogonal projection onto 
$K$. By (4.28) or more precisely its proof, $D_{1}\Phi $ maps onto $K$ 
and hence $Im \pi_{K}D_{1} = K$. Since the finite dimensional space $K$ 
splits, we have $TMet_{AH} = K \oplus  K^{\perp} = Im(\pi_{K}D_{1}\Phi 
) \oplus  Ker(\pi_{K}D_{1}\Phi ),$ so that $Ker(\pi_{K}D_{1}\Phi )$ 
splits. Hence, there is a bounded linear projection $P_{1}$ onto 
$Ker(\pi_{K}D_{1}\Phi)$. The operator $L + \pi_{K}$ is invertible and 
one may now define
$$P(\dot \gamma, \dot h) = (P_{1}\dot \gamma, (L+\pi_{K})^{-1}
(-(D_{1}\Phi)P_{1}(\dot \gamma) + \pi_{K}\dot h)). $$
Then $P$ is the required bounded linear projection.
{\endproof}

  As in Remark 3.2, it is doubtful if Theorem 4.1 remains valid in general 
without the assumption $\pi_{1}(M, \partial M) = 0$. As noted there, in the 
generic situation where $g \in E_{AH}^{k,\beta}$ has no local Killing fields, 
Theorem 4.1 does hold at $g$, at least when $\partial M$ is connected. For 
simplicity, for the rest of this section and throughout \S 5, we assume 
$\pi_{1}(M, \partial M) = 0$. 

\begin{corollary} \label{c 4.2.}
  For any $C^{3,\alpha}$ conformally compact metric $g \in E_{AH}^{k,\beta}$, 
the local space $Z_{AH}^{k,\beta}$ is an infinite dimensional $C^{\infty}$ 
separable Banach manifold. In fact, via the splitting {\rm (4.13)}, 
$Z_{AH}^{k,\beta}$ is a $C^{\infty}$ Banach submanifold of 
$Met^{m,\alpha}(\partial M)\times {\mathbb S}_{2}^{k,\beta}(M)$ and 
as such
\begin{equation} \label{e4.35}
T_{g}Z_{AH}^{k,\beta} = Ker(D\Phi )_{g}. 
\end{equation}
\end{corollary}
{\bf Proof:}
 This is an immediate consequence of the definition (4.20), Theorem 4.1 
and the implicit function theorem in Banach spaces, cf.~[18]. 
$Z_{AH}^{k,\beta}$ is separable since it is a submanifold of 
$Met^{m,\alpha}(\partial M)\times {\mathbb S}_{2}^{k,\beta}(M)$, each of which are 
separable Banach spaces in the topologies defined at the beginning of \S 4. 
{\endproof} 

 Locally, near any given $g_{0} \in E_{AH}^{k,\beta}$, the boundary map taking an 
AH Einstein metric $g$ to its boundary metric $\gamma$ with respect to the 
compactification (4.1) is given simply by projection on the first factor:
\begin{equation} \label{e4.36}
\Pi : E_{AH}^{k,\beta} \rightarrow  Met^{m,\alpha}(\partial M);  \ \ 
\Pi(g) = \Pi (g_{\gamma}+h) = \gamma . 
\end{equation}
Clearly, this map is $C^{\infty}$ smooth.

 The spaces $Met_{AH}^{k,\beta}$ and $E_{AH}^{k,\beta}$ are invariant under 
the action of suitable diffeomorphisms. In \S 5, we will consider larger 
diffeomorphism groups, but for now we restrict to the group ${\mathcal D}_{2} \equiv$ 
Diff$^{m+1,\alpha}(\bar M)$ of $C^{m+1,\alpha}$ diffeomorphisms $\phi$ of 
$\bar M$ such that
\begin{equation} \label{e4.37}
\phi|_{\partial M} = id_{\partial M}, \ \ {\rm and} \ \  
\lim_{\rho_{o}\rightarrow 0}(\frac{\phi^{*}\rho_{0}}{\rho_{0}}) = 1, 
\end{equation}
where $\rho_{0}$ is the fixed defining function. If $g\in 
E_{AH}^{k,\beta}$ and $\widetilde g = \rho_{0}^{2}g$ is the compactification 
as in (4.1), then for $\phi\in{\mathcal D}_{2}$, the compactification of 
$\phi^{*}g$ is given by
$$\widetilde{\phi^{*}g} = \rho_{0}^{2}\phi^{*}g = (\phi^{*}\widetilde g)
(\frac{\rho_{0}}{\phi^{*}\rho_{0}})^{2}. $$
Hence (4.37) implies that $g$ and $\phi^{*}g$ have the same boundary 
metric with respect to $\rho_{0}.$ However, the normal vectors of the 
compactified metrics $\widetilde g$ and $\phi^{*}\widetilde g$ are different 
in general.

 The action of ${\mathcal D}_{2}$ preserves the spaces $Met_{AH}^{k,\beta}$ 
and $E_{AH}^{k,\beta}$. This is because $|D\phi - id|_{\widetilde g}$ 
extends $C^{m,\alpha}$ smoothly up to $\partial M$, and hence 
$|D\phi - id|_{g} = O(e^{- r})$, so that $|\phi^{*}g - g|_{g} = 
O(e^{-2r})$. Note also that since $m \geq k$, for $g$ a $C^{k,\beta}$ metric 
(in a smooth atlas for $M$), and $\phi \in {\mathcal D}_{2}$, $\phi^{*}g \in 
C^{k,\beta}$. 

  Observe that the action of ${\mathcal D}_{2}$ on $Met_{AH}^{k,\beta}$ 
or $E_{AH}^{k,\beta}$ is free, since any isometry $\phi$ of a metric inducing the 
identity on $\partial M$ must itself be the identity; this is most easily seen 
by working in a geodesic compactification $\bar g$. It is also standard that 
the action of ${\mathcal D}_{2}$ on $Met_{AH}^{k,\beta}$ and $E_{AH}^{k,\beta}$ 
is proper. 

  It is well-known however that the action of ${\mathcal D}_{2}$ on $Met_{AH}^{k,\beta}$ 
is not smooth; for a 1-parameter group of diffeomorphisms $\phi_{t}$ with 
$\phi_{0} = id$ and infinitesimal generator $X$, one has 
$\frac{d}{dt}(\phi_{t}^{*}g)|_{t} = \phi_{t}^{*}{\mathcal L}_{X}g$. For 
$g \in Met_{AH}^{k,\beta}$ and $X \in T_{id}{\mathcal D}_{2}$, the form 
${\mathcal L}_{X}g$ is only $C^{k-1,\beta}$ smooth and so not an element of 
$T_{g}Met_{AH}^{k,\beta}$. However, as noted following \eqref{e4.16}, Einstein 
metrics $g$ are $C^{\infty}$ smooth in a smooth atlas for $M$, and in such coordinates, 
${\mathcal L}_{X}g$ is $C^{k,\beta}$, (in fact $C^{m,\alpha}$), smooth. Thus, 
there is no loss-of-derivatives for Einstein metrics. 

   Now it is proved in [8, Prop. I.4.6] that the set of metrics $g \in 
E_{AH}^{k,\beta}$ near a given $g_{0} \in E_{AH}^{k,\beta}$ such that
\begin{equation} \label{e4.38}
\beta_{g_{0}}(g) = 0 
\end{equation}
is a local slice for the action of ${\mathcal D}_{2}$ on $E_{AH}^{k,\beta}$. 
Thus, a neighborhood ${\mathcal U}$ of any given $g_{0}\in E_{AH}^{k,\beta}$ 
is homeomorphic to a product $Z_{AH}^{k,\beta}\times {\mathcal V}$, where 
${\mathcal V}$ is a neighborhood of the identity in ${\mathcal D}_{2}$. 
The homeomorphism is given by 
$$\psi_{0}(g) = (\phi_{0}^{*}g, \phi_{0}),$$
where $\phi_{0}$ is the unique element of ${\mathcal D}_{2}$ such that 
$\phi_{0}^{*}g \in Z_{AH}^{k,\beta}$. To consider the corresponding overlap 
maps, let $g_{0}$ and $g_{1}$ be background metrics in $E_{AH}^{k,\beta}$ 
which are sufficiently close, and let $Z_{i}$ be the space (4.20) determined 
by $g_{i}$. Then $\psi_{0}(g) = (\phi_{0}^{*}g, \phi_{0})$, $\psi_{1}(g) 
= (\phi_{1}^{*}g, \phi_{1})$ with $\phi_{i}^{*}g \in Z_{i}$, and hence the 
overlap map is given by 
$$\psi_{01}(g_{0}, \phi_{0}) = (g_{1}, \phi_{1}) = 
((\phi_{1}^{*}(\phi_{0}^{-1})^{*})g_{0}, (\phi_{1}\circ \phi_{0}^{-1})\phi_{0}),$$
where $g_{i} \in Z_{i}$ and $\phi_{1} = \phi_{1}(g_{0}, \phi_{0})$ is defined 
as the unique solution of the equation $\beta_{g_{1}}(\phi_{1}^{*}(\phi_{0}^{-1})^{*}g_{0}) 
= 0$. By the discussion preceding (4.38), $\phi_{1}$ is differentiable in $g_{0}$ and 
$\phi_{0}$ and in fact is $C^{\infty}$ smooth in these variables. It follows that the 
overlap maps are $C^{\infty}$ and hence the global space $E_{AH}^{k,\beta}$ 
is a $C^{\infty}$ smooth separable Banach manifold, as is the quotient
\begin{equation} \label{e4.39}
{\mathcal E}_{AH}^{(2)} = E_{AH}^{k,\beta}/{\mathcal D}_{2}^{m+1,\alpha}. 
\end{equation}
Two metrics $g_{1}$ and $g_{2}$ in ${\mathcal E}_{AH}^{(2)}$ are equivalent if 
there is a $C^{m+1,\alpha}$ diffeomorphism $\phi$ of weight 2, i.e.~satisfying 
(4.37), such that $\phi^{*}g_{1} = g_{2}$. In particular, $g_{1}$ and 
$g_{2}$ must have the same boundary metric with respect to $\rho_{0}$. 

  When $E_{AH}^{k,\beta}$ is viewed as subset of the product 
$Met^{m,\alpha}(\partial M)\times {\mathbb S}_{2}^{k,\beta}(M)$ via $S$, since 
${\mathcal D}_{2}^{m+1,\alpha}$ acts trivially on the first factor, one has
\begin{equation} \label{e4.40}
{\mathcal E}_{AH}^{(2)} \subset  Met^{m,\alpha}(\partial M)\times 
({\mathbb S}_{2}^{k,\beta}(M)/{\mathcal D}_{2}^{m+1,\alpha}). 
\end{equation}
This inclusion sends $[g] = [g_{\gamma}+h]$ to $(\gamma ,[h])$, and, given 
a fixed $g_{0}\in E_{AH}^{k,\beta}$, a slice representative for $[h]$ is that 
unique $h\in [h]$ satisfying (4.38). Via (4.36), $\Pi$ descends to a smooth map
\begin{equation} \label{e4.41}
\Pi : {\mathcal E}_{AH}^{(2)} \rightarrow  Met^{m,\alpha}(\partial M), \ \ 
\Pi ([g]) = \gamma . 
\end{equation}
We summarize the analysis above in the following:
\begin{proposition} \label{p 4.3.}
  Near any $C^{3,\alpha}$ conformally compact Einstein metric $g$, the space 
${\mathcal E}_{AH}^{(2)}$ is a smooth separable Banach manifold. The map 
$\Pi : {\mathcal E}_{AH}^{(2)} \rightarrow  Met^{m,\alpha}(\partial M)$ is 
a $C^{\infty}$ Fredholm map of index 0, with
\begin{equation} \label{e4.42}
Ker (D\Pi )_{g} = K_{g}, 
\end{equation}
where as in {\rm (4.24)}, $K_{g}$ is the space of $L^{2}$ infinitesimal 
Einstein deformations at g. Consequently, $Im\Pi  \subset  
Met^{m,\alpha}(\partial M)$ is a variety of finite codimension.
\end{proposition}
{\bf Proof:}
 One only needs to verify that $\Pi$ is Fredholm, with kernel given by (4.42). 
By construction, one has
$$Ker D\Pi  = T{\mathcal E}_{AH}^{(2)} \cap  Ker \Pi_{1},$$
where $\Pi_{1}$ is the linear projection on the first factor in the 
splitting (4.40). Since $D\Pi_{1}= Id$ on the first factor,
$$T{\mathcal E}_{AH}^{(2)} \cap  Ker \Pi_{1} = T{\mathcal E}_{AH}^{(2)}\cap 
T({\mathbb S}_{2}^{k,\beta}/{\mathcal D}_{2}^{m+1,\alpha}).$$
This intersection just consists of the classes $[h]$ satisfying (4.38), 
and so by (4.35) and (4.26), 
$$Ker D\Pi  = Ker L, $$
where the kernel is taken in ${\mathbb S}_{2}^{k,\beta}.$ But this is the 
same as the $L^{2}$ kernel $K$ by [8, Prop.I.3.5].

 For the cokernel, one has
$$Im(D\Pi ) = \Pi (T{\mathcal E}_{AH}^{(2)}) = \Pi (KerD\Phi ) = 
Ker(\pi_{K}D_{1}\Phi ), $$
where the second equality is from (4.35) and the last equality follows 
from (4.34) and the discussion following it. Again, as following 
(4.34), $Ker \pi_{K}D_{1}\Phi  = K^{\perp}$ is closed, and has 
codimension $k = dim K$. Hence $\Pi$ is Fredholm of index 0.
{\endproof}
\begin{remark} \label{r 4.4.}
  {\rm This result shows that one has the following dictotomy: either there 
exist no conformally compact Einstein metrics on $M$, or the moduli space of 
such metrics is at least infinite dimensional, with $Im \Pi$ a variety of 
finite codimension in $Met (\partial M)$. 

   If there exist Einstein metrics $g \in E_{AH}$ which are non-degenerate, 
so that $K_{g} = \{0\}$, then $\Pi$ is a local diffeomorphism in a 
neighborhood of $g$. This is the result of Biquard [8], extending earlier 
work of Graham-Lee [15]. In other words, $\Pi$ is an open map on the 
open submanifold $E'_{AH}$ of non-degenerate metrics. 
  
  Note that $T{\mathcal E}_{AH}^{(2)}$ is the space of (essential) 
infinitesimal asymptotically hyperbolic Einstein deformations, (not 
necessarily preserving the boundary metric as is the case with the 
$L^{2}$ kernel $K$). The fact that ${\mathcal E}_{AH}^{(2)}$ is a smooth 
Banach manifold implies that any infinitesimal AH Einstein deformation 
may be integrated to a (local) curve of AH Einstein metrics. Apriori, it 
is not clear if this remains the case when the boundary metric is required 
to be fixed, i.e.~an $L^{2}$ infinitesimal Einstein deformation in $K$ 
might not integrate to a curve of AH Einstein metrics with the same 
boundary metric. }
\end{remark} 

 Observe that all the results above are valid in any dimension.

\medskip

 We complete this section with a discussion of the boundary regularity 
of metrics in $E_{AH}$. The Einstein metrics in $E_{AH}^{k,\beta}$ have 
$C^{2}$ and hence $L^{2,p}$ compactifications. Suppose dim $M = n+1 = 4$. 
Then Theorem 2.3 implies that any $g\in E_{AH}^{k,\beta}$ is 
$C^{m,\alpha}$ conformally compact, for any $m \geq 2$, see the 
discussion following (4.2). Thus
\begin{equation} \label{e4.43}
E_{AH}^{k,\beta} = E_{AH}^{m,\alpha}, 
\end{equation}
and $E_{AH}^{m,\alpha}$ is the space of AH Einstein metrics on $M$ 
which are $C^{m,\alpha}$ conformally compact with respect to the defining 
function $\rho_{0}$ as in (4.1). The space $E_{AH}^{m,\alpha}$ 
is a smooth separable Banach manifold, and boundary regularity implies that 
the topology on $E_{AH}^{m,\alpha}$ defined by (4.16) is equivalent to the 
${\mathbb S}_{2}^{m,\mu}(\bar M)$ topology on the compact manifold $\bar M$, 
for a fixed $\mu < \alpha$. This corresponds to the definition in the Introduction. 
 With this understood, one has the following version of Theorem 1.2:

\begin{proposition} \label{p 4.5.}
  If dim $M = 4$ and $m \geq 3$, then $E_{AH}^{m,\alpha}$ is the space of 
$C^{m,\alpha}$ conformally compact Einstein metrics on $M$. If 
$\pi_{1}(M, \partial M) = 0$, then $E_{AH}^{m,\alpha}$ is a smooth 
separable Banach manifold and the map $\Pi$ is a $C^{\infty}$ map
\begin{equation} \label{e4.44}
\Pi : E_{AH}^{m,\alpha} \rightarrow  Met^{m,\alpha}(\partial M). 
\end{equation}
\end{proposition}
{\endproof}

 Of course, Proposition 4.5 also holds on the quotient ${\mathcal E}_{AH}^{(2)} 
= {\mathcal E}_{AH}^{(2),m,\alpha}.$ An analogous but somewhat weaker result 
holds in all higher dimensions $n+1 > 4$; in fact there are two versions in 
higher dimensions, although neither version is quite as strong as Proposition 
4.5, cf.~Theorems 5.5 and 5.6 for further details.

\section{The Spaces ${\mathcal E}^{m,\alpha}$, Diffeomorphisms and 
Stability.}
\setcounter{equation}{0}

 In \S 4, the defining function $\rho_{0}$ was fixed, thus giving a 
fixed boundary metric $\gamma $ for an AH Einstein metric $g$ on $M$. 
In this section, we consider the situation where $\rho$ varies over all 
smooth defining functions, and the corresponding variation of the 
boundary metrics. This is closely related to the action of 
diffeomorphisms on $\bar M$. These issues are discussed in \S 5.1, 
together with the proof of Theorem 1.2 and its versions in higher 
dimensions. In \S 5.2, we prove that the spaces ${\mathcal 
E}_{AH}^{m,\alpha}$ are all diffeomorphic and stable in a natural 
sense. 

{\bf \S 5.1.}
 Let ${\mathcal D}_{1} = {\mathcal D}_{1}^{m+1,\alpha}(\bar M)$ be the group 
of orientation preserving $C^{m+1,\alpha}$ diffeomorphisms of $\bar M$ 
which restrict to the identity map on $\partial M$. Recall from (4.37) 
that ${\mathcal D}_{2} \subset {\mathcal D}_{1}$ is the subgroup of 
diffeomorphisms $\phi$ satisfying $\lim_{\rho_{0}\rightarrow 0}
(\phi^{*}\rho_{0}/\rho_{0}) = 1$. It is easily seen that ${\mathcal D}_{2}$ 
is a normal subgroup of ${\mathcal D}_{1}$. With respect to $\rho_{0}$, one 
has the splitting $TM|_{\partial M} \cong  T(\partial M)\oplus{\mathbb R}$, 
where the ${\mathbb R}$ factor is identified with the span of 
$\partial /\partial\rho_{0}$. The groups ${\mathcal D}_{1}$ and ${\mathcal D}_{2}$ 
act on $T(\partial M)\oplus{\mathbb R}$ by the map $\phi \rightarrow 
D\phi|_{\partial M}$, and so induce subgroups of $Hom(TM|_{\partial M}, 
TM|_{\partial M})$. Since ${\mathcal D}_{2} \subset {\mathcal D}_{1}$ is 
defined solely by a $1^{\rm st}$ order condition at $\partial M$, the 
quotient group ${\mathcal D}_{1}/{\mathcal D}_{2}$ is isomorphic to the 
corresponding quotient group in $Hom(TM|_{\partial M}, TM|_{\partial M})$.

\begin{lemma} \label{l 5.1.}
  The quotient group ${\mathcal D}_{1}/{\mathcal D}_{2}$ is naturally 
isomorphic to the group of $C^{m,\alpha}$ positive functions $\lambda$ 
on $\partial M$.
\end{lemma}
{\bf Proof:}
 With respect to the splitting $TM|_{\partial M} \cong  T(\partial M)
\oplus{\mathbb R}$, the linear map $D\phi|_{\partial M}$, for 
$\phi\in{\mathcal D}_{1}$, has the form
$$\left( \begin{array}{cc} 
       1 & * \\
       0 & \lambda 
\end{array} \right)$$
where $\lambda  = \lim_{\rho_{0}\rightarrow 0}(\phi^{*}\rho_{0}/\rho_{0})$. 
For $\phi\in{\mathcal D}_{2}$, $D\phi$ is the same, except that the entry 
$\lambda$ is $1$. It follows that the quotient group is identified with the 
multiplicative group of functions $\lambda : {\mathbb R}  \rightarrow  
{\mathbb R}$, acting in the $\partial /\partial\rho_{0}$ direction. Since 
$D\phi$ is non-singular, $\lambda$ cannot vanish and hence $\lambda  > 0$.
{\endproof} 

  As in \S 4, let ${\mathcal E}_{AH}^{(2)} = E_{AH}/{\mathcal D}_{2}$ be the 
space of isometry classes of AH Einstein metrics, among diffeomorphisms 
in ${\mathcal D}_{2}$, and similarly, let ${\mathcal E}_{AH}^{(1)} = 
E_{AH}/{\mathcal D}_{1}$; here $E_{AH} = E_{AH}^{k,\beta}$, as in (4.16), 
(or (4.43)). There is a natural projection map 
${\mathcal E}_{AH}^{(2)}\rightarrow {\mathcal E}_{AH}^{(1)}$ with fiber 
${\mathcal D}_{1}/{\mathcal D}_{2}$. As in \S 4, ${\mathcal D}_{1}$ acts freely on 
$E_{AH}$, with local Bianchi slice as in (4.38) so that as following (4.38), 
${\mathcal E}_{AH}^{(1)}$ is a $C^{\infty}$ separable Banach manifold.

 Next, let ${\mathcal C}  = {\mathcal C}^{m,\alpha}$ be the space of conformal 
classes of $C^{m,\alpha}$ metrics on $\partial M$. Again, ${\mathcal C}$ 
has the structure of an infinite dimensional Banach manifold, with 
tangent spaces given by the space of trace-free symmetric bilinear 
forms. There is a natural projection map $Met^{m,\alpha}(\partial M) 
\rightarrow  {\mathcal C}$, with fiber the space of $C^{m,\alpha}$ 
conformally equivalent metrics on $\partial M$.
\begin{proposition} \label{p 5.2.}
  The boundary map $\Pi$ descends to a $C^{\infty}$ boundary map on 
the base spaces, i.e.
\begin{equation} \label{e5.1}
\Pi : {\mathcal E}_{AH}^{(1)} \rightarrow  {\mathcal C} . 
\end{equation}
This map $\Pi$ is Fredholm, of index 0, with Ker $D\Pi = K$, as in 
{\rm (4.42)}.
\end{proposition}
{\bf Proof:}
 Let $g_{1}$ and $g_{2}$ be AH Einstein metrics on $M$ with 
$\phi^{*}g_{2} = g_{1}$, for $\phi\in{\mathcal D}_{1},$ and set $\lambda  = 
\lim_{\rho_{0}\rightarrow 0}(\phi^{*}\rho_{0}/\rho_{0})$. Let $\bar g_{i}$ 
be the compactification of $g_{i}$, $i = 1,2$, with respect to $\rho_{0}$, 
as in (4.1), and let $\gamma_{i}$ be the induced boundary metrics. If 
$\hat g_{2}$ is the $\rho_{0}$-compactification of $\phi^{*}g_{2}$, 
then one has
$$\hat g_{2} = \rho_{0}^{2}\phi^{*}(\rho_{0}^{-2})\phi^{*}(\rho_{0}^{2}g_{2}) = 
(\frac{\rho_{0}}{\phi^{*}(\rho_{0})})^{2}\phi^{*}(\rho_{0}^{2}g_{2}). $$
Hence, the boundary metric $\hat \gamma_{2}$ of $\hat g$, which must 
equal $\gamma_{1},$ is given by
$$\gamma_{1} = \hat \gamma_{2} = \lambda^{-2}\phi^{*}\gamma_{2}. $$
Since $\phi  = id$ on $\partial M$, it follows that $\gamma_{2} = 
\lambda^{2}\gamma_{1}$, so that the boundary metrics are conformal. It 
follows that the boundary map $\Pi$ in (4.41) descends to the map $\Pi$ 
in (5.1) and is smooth.

  Further, observe that Lemma 5.1 shows that the converse of the proof 
above also holds, i.e.~if $\gamma_{1}$ and $\gamma_{2}$ are conformally 
equivalent metrics in $Met^{m,\alpha}(\partial M)$, so that $\gamma_{2} 
= \lambda^{2}\gamma_{1}$, then there is a diffeomorphism $\phi\in{\mathcal D}_{1}$ 
such that $\phi^{*}g_{2} = g_{1}$, where $g_{i}$ are any AH Einstein metrics 
on $M$ with boundary metrics $\gamma_{i}$ with respect to the 
$\rho_{0}$-compactification. Hence, $\Pi$ maps the fibers 
${\mathcal D}_{1}/{\mathcal D}_{2}$ diffeomorphically onto the fibers of 
$Met^{m,\alpha}(\partial M) \rightarrow  {\mathcal C}$. 

 The proof that $\Pi$ is Fredholm of index 0, with $Ker D\Pi = K$, is 
thus exactly the same as in Proposition 4.3.
{\endproof}
\begin{remark} \label{r 5.3.}
  {\rm Recall that the map $\Pi: {\mathcal E}_{AH}^{(2)} \rightarrow  
Met(\partial M)$ in (4.41) depends on a choice of the defining function 
$\rho_{0}$ from (4.1). The reduced map $\Pi$ in (5.1) is now 
independent of the choice of $\rho_{0}$. To see this, let $\rho_{1}$ be 
any other defining function, so that $\rho_{1} = \lambda\cdot \rho_{0}$, 
for some function $\lambda  > 0$ on $M$. Let
$$\bar g = \rho_{0}^{2}g, \ \ {\rm and} \ \ \widetilde g = \rho_{1}^{2}g $$
be compactifications of $g$ with respect to $\rho_{0}$ and $\rho_{1}$. The 
boundary metrics are related by $\widetilde \gamma = \lambda^{2}\gamma$, 
where $\lambda  = \lim_{\rho\rightarrow 0}(\rho_{1}/\rho_{0})$. As in 
the proof of Proposition 5.2. there is a diffeomorphism $\phi\in{\mathcal D}_{1}$ 
satisfying, (along integral curves of $\partial /\partial\rho_{0}$), 
$d\phi (\rho_{0})/d\rho_{0} = \lambda$ at $\partial M$. Hence
$$\widetilde {(\phi^{*}g)} = \rho_{1}^{2}\phi^{*}g = 
\lambda^{2}\rho_{0}^{2}\phi^{*}g, $$
while
$$\phi^{*}\widetilde g = \phi^{*}(\rho_{0}^{2})\phi^{*}g = 
\lambda^{2}\rho_{0}^{2}\phi^{*}g, $$
near $\partial M$. Thus, the $\rho_{1}$ compactification of $\phi^{*}g$ 
is the same as the $\rho_{0}$ compactification of $g$, pulled back by 
$\phi$. }
\end{remark}

 Theorem 1.2 is now essentially an immediate consequence of the work 
above and in \S 4.

{\bf Proof of Theorem 1.2.}

 The discussion following Lemma 5.1 shows that ${\mathcal E}_{AH}^{(1)} = 
{\mathcal E}_{AH}^{(1),k,\beta}$ and ${\mathcal C} = {\mathcal C}^{m,\alpha}$, 
are $C^{\infty}$ smooth separable Banach manifolds and by Proposition 5.2, 
$\Pi$ is a $C^{\infty}$ Fredholm map of index 0. The boundary regularity 
result in Proposition 4.5, cf.~(4.43), shows that by setting $(k,\beta) = (m,\alpha)$, 
${\mathcal E}_{AH}^{(1)}$ is the space of AH Einstein metrics on $M$ which admit a 
$C^{m,\alpha}$ compactification, with topology that of ${\mathbb S}_{2}^{m,\mu}(\bar M)$, 
$\mu < \alpha$.  
{\endproof}
\begin{remark} \label{r 5.4.}
  {\rm For certain purposes, it is useful to consider quotients by larger 
diffeomorphism groups, and we discuss this briefly here. Thus, let 
${\mathcal D}_{0}$ be the group of $(C^{m+1,\alpha})$ diffeomorphisms of 
$\bar M$ such that the induced diffeomorphism on $\partial M$ is 
isotopic to the identity. Again, the group ${\mathcal D}_{1} \subset 
{\mathcal D}_{0}$ is a normal subgroup, and one may form 
\begin{equation} \label{e5.2}
{\mathcal E}_{AH}^{(0)} = E_{AH}/{\mathcal D}_{0} = 
{\mathcal E}_{AH}^{(1)}/({\mathcal D}_{0}/{\mathcal D}_{1}). 
\end{equation}
Similarly, let ${\mathcal T}$ denote the quotient space ${\mathcal T} = 
{\mathcal C}/{\mathcal D}_{0}$. This is the space of {\it marked} conformal 
structures on $\partial M$, analogous to the Teichm\"uller space of 
conformal structures on surfaces. The group ${\mathcal D}_{0}$ however does 
not act freely on ${\mathcal C}$. Elements in ${\mathcal C}$ having a non-trivial 
isotropy group ${\mathcal D}_{0}[\gamma]$ are the classes $[\gamma]$ which have 
a non-trivial group of conformal diffeomorphisms, i.e.~${\mathcal D}_{0}[\gamma]$ 
consists of diffeomorphisms $\phi\in{\mathcal D}_{0}$ such that
$$\phi^{*}\gamma  = \lambda^{2}\cdot \gamma , $$
for some positive function $\lambda$ on $\partial M$. A well-known 
theorem of Obata [24] implies that the isotropy group ${\mathcal D}_{0}[\gamma]$ 
of $[\gamma]$ is always compact, with the single exception of $(\partial M, [\gamma]) 
= (S^{n-1}, [\gamma_{0}])$, where $\gamma_{0}$ is the round metric on $S^{n-1}$. 

 Similarly, the elements $g$ of $E_{AH}$ which have non-trivial 
isotropy groups ${\mathcal D}_{0}(g)$ in ${\mathcal D}_{0}$ are AH Einstein 
metrics which have a non-trivial group of isometries. Such isometries 
$\phi\in{\mathcal D}_{0}$ induce a diffeomorphism $\phi $ of $\partial M$, 
which is a conformal isometry of the conformal infinity $[\gamma]$ of 
$g$. It follows that the boundary map $\Pi$ in (5.1) descends further to 
a boundary map
\begin{equation} \label{e5.3}
\Pi : {\mathcal E}_{AH}^{(0)} \rightarrow  {\mathcal T} . 
\end{equation}
At any class $[g]$ where ${\mathcal D}_{0}[g] = id$, the quotient space 
${\mathcal E}_{AH}^{(0)}$ is a smooth infinite dimensional Banach manifold, 
and similarly for ${\mathcal T}$. At those classes $[g]$ or $[\gamma]$ 
where ${\mathcal D}_{0}[g]$ or ${\mathcal D}_{0}[\gamma]$ is compact, the 
quotients ${\mathcal E}_{AH}^{(0)}$ and ${\mathcal T}$ are smooth orbifolds, 
and $\Pi$ is an orbifold smooth map. Only at the exceptional class 
$(B^{n}, g_{-1})$ of the Poincar\'e metric on the ball is the quotient 
${\mathcal T}$ not well-behaved, and possibly non-Hausdorff.

 Finally, one may carry out the same quotient construction with respect 
to the full group ${\mathcal D}$  = Diff$(\bar M)$ of diffeomorphisms of 
$\bar M$ mapping $\partial M$ to itself, so that ${\mathcal E}_{AH} = 
E_{AH}/{\mathcal D}$, while ${\mathcal T}$ is replaced by the moduli space of 
conformal structures ${\mathcal M}  = {\mathcal T} /\Gamma$, where $\Gamma  = 
{\mathcal D}(\partial M)/{\mathcal D}_{0}(\partial M)$ is the subgroup of the 
mapping class group of $\partial M$ consisting of diffeomorphisms of 
$\partial M$ which extend to diffeomorphisms of $M$. }
\end{remark}

 Next we discuss two versions of Theorem 1.2 in higher dimensions. Let 
$M$ be an $(n+1)$ dimensional manifold with boundary, $n > 3$. When $n$ is 
even, the Fefferman-Graham expansion (3.3) in general has $\log$ terms 
appearing at order $n$, i.e.~of the form $t^{n}\log t$, and at higher 
order as well. Thus, one cannot expect a smooth boundary regularity result 
when dim $M$ is odd. On the other hand, a result of Lee [20] gives boundary 
regularity below order $n$.

 To describe the first version of Theorem 1.2, let $E_{AH}^{2}$ be the 
space of AH Einstein metrics which are $C^{2}$ conformally compact, 
with respect to a smooth defining function $\rho_{0}$, as in (4.16). 
Suppose the boundary metric $\gamma \in  C^{m,\alpha}$. Then Lee's 
result [20] states that any $g\in  E_{AH}^{2}$ is $C^{k,\mu}$ conformally 
compact, where $k+\mu  = \min(m+\alpha , n-1+\beta)$, for any $\beta \in (0,1)$. 

 Combining this result with the results above in \S 5, and with 
Proposition 4.3 and the discussion preceding Proposition 4.5, gives the following:

\begin{theorem} \label{t 5.5.}
  Let $M$ be a compact, oriented $(n+1)$-manifold with boundary $\partial M$, 
$n > 3$, with $\pi_{1}(M, \partial M) = 0$. If, for a given $(m, \alpha)$, with 
$3 \leq  m \leq n-1$, $\hat {\mathcal E}_{AH} = \hat {\mathcal E}_{AH}^{m,\alpha}$ 
is non-empty, then $\hat {\mathcal E}_{AH}$ is a smooth infinite dimensional Banach 
manifold. Further, the boundary map 
\begin{equation} \label{e5.4}
\Pi : \hat {\mathcal E}_{AH} \rightarrow  {\mathcal C}  = {\mathcal C}^{m,\alpha} 
\end{equation}
is a $C^{\infty}$ smooth Fredholm map of index 0.
\end{theorem}{\endproof} 
 
 Thus, the statement of Theorem 5.5 is equivalent to that of Theorem 
1.2, provided $m \leq  n-1$. For the second version of Theorem 1.2, a 
result of Chru\'sciel et al.~[10] gives an optimal boundary regularity 
result for $C^{\infty}$ boundary metrics $\gamma$. Thus, if $g$ is an AH 
Einstein metric with a $C^{2}$ conformal compactification to a $C^{\infty}$ 
boundary metric $\gamma$, then if $n$ is odd, $g$ is $C^{\infty}$ conformally 
compact. If $n$ is even, $g$ is $C^{\infty}$ polyhomogeneous, i.e.~$g$ has a 
compactification $\bar g$ which is a smooth function of $(t, t^{n}\log t, y)$, 
where $y\in\partial M$. In either case even/odd, let $\widetilde {\mathcal E}_{AH}$ 
be the space of such metrics, and ${\mathcal C}^{\infty}$ the space of $C^{\infty}$ 
conformal classes.

  The same proof as Theorem 5.5 gives:
\begin{theorem} \label{t 5.6.}
  Let $M$ be a compact, oriented $(n+1)$-manifold with boundary 
$\partial M$ with $\pi_{1}(M, \partial M) = 0$. If $\widetilde {\mathcal E}_{AH}$ 
is non-empty, then $\widetilde {\mathcal E}_{AH}$ is a smooth infinite dimensional 
Frechet manifold. Further, the boundary map 
\begin{equation} \label{e5.5}
\Pi : \widetilde {\mathcal E}_{AH} \rightarrow  {\mathcal C}  = {\mathcal C}^{\infty} 
\end{equation}
is a $C^{\infty}$ smooth Fredholm map of index 0.
\end{theorem}
{\endproof} 

 {\bf \S 5.2.}
 In this section, we compare the structure of the spaces ${\mathcal E}_{AH}^{m,\alpha}$ 
over varying $m$, $\alpha$. In dimension 4, the spaces ${\mathcal E}_{AH}^{m,\alpha}$ 
are defined as in (4.43), while in dimensions greater than 4, 
${\mathcal E}_{AH}^{m,\alpha}$ is defined as  preceding (5.4) with $3 \leq m \leq n-1$. 
(The spaces $\widetilde {\mathcal E}_{AH}$ are only defined for $C^{\infty}$ boundary 
data).

 Clearly, one has inclusions
\begin{equation} \label{e5.6}
{\mathcal E}_{AH}^{\omega} \subset {\mathcal E}_{AH}^{\infty} \subset  
{\mathcal E}_{AH}^{m' ,\alpha'}\subset {\mathcal E}_{AH}^{m,\alpha}, 
\end{equation}
for any $(m',\alpha')$ with $m'+\alpha' > m+\alpha$. Here we recall that 
the topology on ${\mathcal E}_{AH}^{m,\alpha}$ is that induced by the 
$C^{m,\mu}$ topology on ${\mathbb S}_{2}^{m,\alpha}(\bar M)$, for a fixed 
$\mu < \alpha$, cf.~the discussion preceding Proposition 4.4. 

 The inclusions (5.6) correspond formally to the much simpler inclusions of the 
conformal classes ${\mathcal C}^{\omega} \subset {\mathcal C}^{\infty} \subset 
 {\mathcal C}^{m',\alpha'}\subset {\mathcal C}^{m,\alpha}$ of conformal 
classes of metrics on $\partial M$. It is essentially clear that the 
spaces ${\mathcal C}^{m,\alpha}$ are diffeomorphic, for all $(m,\alpha)$, 
including $m = \infty $ or $m = \omega$. Further, each ${\mathcal C}^{m',\alpha'}$ 
is dense in ${\mathcal C}^{m,\alpha}$. 

  For later purposes, it is worthwhile to verify these claims explicitly. With 
respect to a fixed real-analytic atlas for $\partial M$, metrics in 
$Met^{m,\alpha}(\partial M)$ are given by a collection of 
$C^{m,\alpha}$ functions $g_{ij}: U \rightarrow  {\mathbb R}$, where $U$ 
is an open set in ${\mathbb R}^{n}$. Hence the topology on 
$Met^{m,\alpha}(\partial M)$ is determined by the standard 
$C^{m,\alpha}$ topology on $C^{m,\alpha}(U, {\mathbb R})$. These local 
spaces are all diffeomorphic in a natural sense, as $(m, \alpha )$ vary 
and induce diffeomorphisms of the global spaces 
$Met^{m,\alpha}(\partial M).$ This argument also holds when passing to 
the associated spaces ${\mathcal C}^{m,\alpha}$ of conformal classes. The 
fact that ${\mathcal C}^{m,\alpha}$ is dense in ${\mathcal C}^{m',\alpha'}$ 
also follows from the fact that the local spaces $C^{m',\alpha'}(U, {\mathbb R})$ 
are dense in $C^{m,\alpha}(U, {\mathbb R})$.

\begin{theorem} \label{t 5.7.}
  For any $(m,\alpha)$, with $m \geq 3$ and including $m = \infty$ and 
$m = \omega$, the spaces ${\mathcal E}_{AH}^{m,\alpha}$ are all diffeomorphic. 
Further ${\mathcal E}_{AH}^{\omega}$, and hence ${\mathcal E}_{AH}^{m',\alpha'}$, 
is dense in ${\mathcal E}_{AH}^{m,\alpha}$ so that if 
$\overline{{\mathcal E}_{AH}^{\omega}}$ denotes the completion of 
${\mathcal E}_{AH}^{\omega}$ in ${\mathcal E}_{AH}^{m,\alpha}$, then
\begin{equation} \label{e5.7}
\overline{{\mathcal E}_{AH}^{\omega}} = {\mathcal E}_{AH}^{m,\alpha}. 
\end{equation}
\end{theorem}
{\bf Proof:}
 It suffices to work with the spaces ${\mathcal E}_{AH}^{(2),m,\alpha}$ in 
(4.39) and $Met^{m,\alpha}(\partial M)$, using a fixed $C^{\omega}$ 
defining function $\rho_{0}$ as in \S 4. In the following, we will drop 
the superscript (2) from the notation. Suppose first that 
$g\in{\mathcal E}_{AH}^{m,\alpha}$ is a regular point of $\Pi$, so that 
$D\Pi_{g}$ is an isomorphism. The inverse function theorem implies that 
there are neighborhoods ${\mathcal U}$ of $g$ in ${\mathcal E}_{AH}^{m,\alpha}$ and 
${\mathcal V}$ of $\gamma = \Pi (g)$ in $Met^{m,\alpha}(\partial M)$ such 
that $\Pi : {\mathcal U} \rightarrow {\mathcal V}$ is a diffeomorphism. 
Since ${\mathcal V}$ is an open set in a Banach space, $\Pi|_{{\mathcal U}}$ 
is a chart for ${\mathcal E}_{AH}^{m,\alpha}$. It follows that 
${\mathcal V}^{m',\alpha'} = Met^{m',\alpha'}(\partial M)\cap{\mathcal V} \subset 
Im \Pi$ and by boundary regularity that 
$\Pi^{-1}({\mathcal V}^{m' ,\alpha'}) \cap {\mathcal U} = {\mathcal U}^{m',\alpha'}$ is 
an open set in ${\mathcal E}_{AH}^{m' ,\alpha'}$. Hence $\Pi$ induces a chart for 
${\mathcal E}_{AH}^{m' ,\alpha'}$ and so ${\mathcal E}_{AH}^{m' ,\alpha'}$ is locally 
diffeomorphic to ${\mathcal E}_{AH}^{m,\alpha}$. 

 Next suppose that $g$ is a singular point of $\Pi$, and let 
$K = Ker D\Pi_{g} \subset  T_{g}Met_{AH}^{m,\alpha}(M)$, with 
$H = (Coker D\Pi_{g})^{\perp} \subset  T_{\Pi(g)}Met_{AH}^{m,\alpha}(\partial M)$, 
where the orthogonal complement is taken with respect to the $L^{2}$ inner 
product. By the implicit function theorem, i.e.~Theorem 4.1, a neighborhood 
${\mathcal U}$ of $g$ in ${\mathcal E}_{AH}^{m,\alpha}$ may be written as a graph 
over a domain ${\mathcal V}$ in $K\oplus H$. This gives a local chart on ${\mathcal U}$ 
and for the same reasons as above, one thus obtains a local chart structure for 
the open set ${\mathcal U}^{m',\alpha'} = {\mathcal E}_{AH}^{m',\alpha'}\cap{\mathcal U} 
\subset {\mathcal E}_{AH}^{m',\alpha'}$.

 These local chart structures patch together to give the spaces 
${\mathcal E}_{AH}^{m,\alpha}$ the Banach manifold structure. Since the 
local charts for ${\mathcal E}_{AH}^{m',\alpha'}$ are just those obtained 
by restricting the charts of ${\mathcal E}_{AH}^{m,\alpha}$ to subdomains, 
it follows that the spaces ${\mathcal E}_{AH}^{m,\alpha}$ are all 
diffeomorphic. Similarly, (5.7) follows from the density of the 
corresponding local charts, i.e.~the density of $C^{\omega}(U, {\mathbb R})$ 
in $C^{m,\alpha}(U, {\mathbb R})$.

{\endproof}

\bibliographystyle{plain}

\smallskip

\begin{center}
February 2008/October 2009
\end{center}

\medskip

\noindent
\address{Department of Mathematics\\
S.U.N.Y. at Stony Brook\\
Stony Brook, NY 11794-3651}\\
\email{anderson@math.sunysb.edu}

\end{document}